\documentclass[12pt]{amsart}
\usepackage{amssymb}

\newtheorem{theorem}{Theorem}
\newtheorem{example}[theorem]{Example}
\newtheorem{proposition}[theorem]{Proposition}
\newtheorem{corollary}[theorem]{Corollary}
\newtheorem{remark}[theorem]{Remark}
\newtheorem{definition}[theorem]{Definition}
\newtheorem{lemma}[theorem]{Lemma}

\newcommand{\E}{{\mathcal E}}

\newcommand{\K}{{\mathcal K}}
\newcommand{\B}{{\mathfrak B}}
\newcommand{\X}{{\mathfrak W}}
\newcommand{\Z}{{\mathbb Z}}
\newcommand{\Q}{{\mathbb Q}}
\newcommand{\R}{{\mathbb R}}
\newcommand{\C}{{\mathbb C}}

\newcommand{\RP}{{\mathbb RP}}
\newcommand{\CP}{{\mathbb CP}}

\begin{document}

\title{Cohomology of spaces of complex knots}

\author{V.A.~Vassiliev}\footnote{Weizmann Institute of Science, Rehovot, Israel.  \\ Email victor.vasilyev@weizmann.ac.il 
 }


\begin{abstract}
We develop a technique for calculating the cohomology groups of spaces of complex parametric knots in \ $\C^k$, \ $k \geq 3$, \ and obtain these groups of low dimensions.
\end{abstract}

\maketitle

\section{Introduction}

V.I.~Arnold asked (see \cite{A}, Problem 1998-10) how to ``complexify'' the theory of Vassiliev invariants. We propose here an answer to this question by a straightforward extension of constructions from \cite{ks}, \cite{twocon} to the case of complex parametric curves in affine spaces. Namely, we describe a method of calculating the cohomology groups of spaces of such curves that do have no cusps or self-intersections. For another view of the Arnold's problem (in the context of the ``complexification'' of the Gauss linking number), see \cite{KhR}.
\medskip

Let \ $d$ \ be a natural number, and \ $P(d,k)$ \ be the space of maps \ $\C^1 \to \C^k$, \ $k \geq 3$, \ defined by arbitrary \ $k$ \ polynomials in variable \ $z$ \ of the form
\begin{equation}
\label{param}
z^d + \lambda_1 z^{d-1} + \dots + \lambda_{d-1} z + \lambda_d , 
\end{equation} 
$\lambda_j \in \C^1$.
 This space is obviously diffeomorphic to $\R^{2k d}$.

A map \ $f: \C^1 \to \C^k$ \ of this type is called a {\em complex knot} of degree \ $d$ \ if it is a smooth embedding, i.e. has no self-intersections and no points of vanishing derivative. The set of maps \ $f \in P(d,k)$ \ that are not complex knots is called the {\em discriminant} and is denoted by \ $\Sigma$.

Below we start computing the cohomology groups of the spaces \ $P(d,k) \setminus \Sigma$ \ for sufficiently large \ $d$ \ using a simplicial resolution of the discriminant space. 

This method also allows us to prove the following stabilization theorem.

\begin{theorem}
\label{stabili} For any natural numbers \ $s$ \ and \ $k \geq 3$, \ there is a number \ $d_k(s)$ \ such that there are natural isomorphisms
 \begin{equation}
\label{stabiso}
H^i(P(d',k) \setminus \Sigma) \simeq H^i(P(d_k(s),k) \setminus \Sigma)
\end{equation}
for all \ $d'>d_k(s)$ \ and any \ $i \leq s$.
\end{theorem}

Thus, for any \ $k \geq 3$ \ the stable cohomology ring of spaces \ $P(d,k) \setminus \Sigma$ \ with \ $d \to \infty$ \ is well defined.
It is natural to consider it as the cohomology ring of the space of complex knots in \ $\C^k$. \ In contrast to the real case, all its finite-dimensional elements are definitely of finite type (while the analogous statement for real knots is an uncertain conjecture).

\begin{table}
\caption{}
\label{tt0}
\begin{tabular}{|c||c|c|c|c|c|c|c|}
\hline
$i$ & $0$ & $2k-5$ & $4k-9$ & $4k-7$ & $6k-14$ & $6k-13$ & $6k-12$ \\
\hline
$H^i$ & $\Z$ & $\Z$ & $\Z_2$ & $\Z_3$ & $\Z_2$ & $\Z_3$ & $\Z_3\oplus \Z$ \\
\hline
\end{tabular}\ 
\end{table}

\begin{theorem}
\label{mthm} If \ $d$ \ is sufficiently large and \ $k>3$, \ then

\begin{enumerate}
\item 
all non-trivial groups \ $H^i(P(d,k) \setminus \Sigma, \Z)$ \ with \ $ i \leq 6k-12$ \ are as shown in Table \ref{tt0}.
 
\item rational homology groups \ $H^i(P(d,k) \setminus \Sigma, \Q)$ \ with \
$i \leq 8k-17$ \ are non-trivial only for \ $i=0,$ $2k-5, $ $6k-12,$ $6k-9$, \ and \ $8k-17$; \ for \ $k>4$ \ these groups with \ $i = 0, 2k-5, 6k-12$ \ and \ $6k-9$ \ are one-dimensional, while if \ $k=4$ \ then \ $6k-9=8k-17$ \ and the corresponding group \ $H^{15}(P(d,4) \setminus \Sigma, \Q)$ \ is at least two-dimensional.
\end{enumerate}

\noindent
If \ $d$ \ is sufficiently large, \ then 

\begin{enumerate}
\item
the group \ $H^i(P(d,3) \setminus \Sigma, \Z)$ \ is equal to \ $\Z$ \ for \ $i = 0 $ \ or \ $1$, \ to \ $\Z_2$ \ for \ $i = 3$ \ or \ $4$, \ is trivial for \ $i=2$, and consists of $18$ elements for \ $i=5$;  

\item
all non-trivial groups \ $H^i(P(d,3) \setminus \Sigma, \Q)$ \ with \ $i \leq 7$ \ are \ $H^0$, $H^1$, $ H^6$ \ and \ $H^7$; \ the first three of them are one-dimensional.
\end{enumerate}
\end{theorem}

\subsection{Notation} For a topological space $X$, \ $B(X,k)$ \ is its \ {\em $k$-th configuration space,} i.e., the space of subsets of cardinality \ $k$ \ in \ $X$ \ with a natural topology. \ $\pm \Z$ \ is the {\em sign local system} of groups on the space \ $B(X,k)$: \ it is locally isomorphic to \ $\Z$, \ but loops in \ $B(X,k)$ \ act on its fibers by multiplication by \ $\pm 1$ \ depending on the parity of corresponding permutations of \ $k$ \ points. \
 $\bar H_*(X)$ \ denotes the Borel--Moore homology group of the topological space \ $X$, \ that is, the homology group of the complex of locally finite singular chains in \ $X$.

\subsection{Work plan}
In \S \ref{sr} we introduce the main tool of the work, a spectral sequence arising from the simplicial resolution of the discriminant variety. Its construction almost repeats that which is systematically used in \cite{book}. Theorem \ref{stabili} is proved in \S \ref{prostab} using this spectral sequence. In \S \ref{prelim} a crucial technical tool of practical calculations is described: an additional filtration on the resolution spaces, which simplifies the calculation of the first page \ $E_1$ \ of the main spectral sequence.

In \S\S \ref{frst}--\ref{sec3t} we apply these techniques to find three first columns of this page, which provide almost all cohomology classes mentioned in Theorem \ref{mthm}. In \S \ref{estim} we study some configuration spaces and local systems involved in the calculations of this kind. In particular, we prove that all subsequent columns almost do not contribute to cohomology groups of low dimensions mentioned in Theorem \ref{mthm} (except for a class of \ $H^{8k-17}(P(d,k) \setminus \Sigma, \Q)$ \ coming from the fourth column). 

\section{Simplicial resolution and main spectral sequence}
\label{sr}

\subsection{Systems of elementary conditions and the first reductions}

Denote by \ ${\X}$ \ the space \ $\mbox{Sym}^2(\C^1)$ \ of unordered pairs of points \ $\{\alpha, \beta\}$ \ in \ $\C^1$. By Vieta theorem, it can be identified with the 
space of polynomials of the form \ $t^2 + u t + v$ \ in the variable \ $t$ \ with complex coefficients, in particular, is homeomorphic to \ $\C^2$.

\begin{definition} \rm
 For any point \ $\chi = \{\alpha, \beta\} \in \X$, \ the corresponding {\em elementary condition} on the maps \ $f: \C^1 \to \C^k$ \ is the condition \ $f(\alpha) = f(\beta)$ \ if \ $\alpha \neq \beta$ \ or \ $f'(\alpha)=0$ \ if \ $\alpha=\beta$. \ The space \ $\X$ \ is called the {\em space of elementary conditions}.

$\divideontimes \subset \Sigma$ \ is the set of all maps of class \ $P(d,k)$ \ satisfying infinitely many elementary conditions.

For any natural \ $d$, \ the \ {\em $d$-rank} of a finite system of elementary conditions \ $\chi_j = \{\alpha_j, \beta_j\}$ \ is the complex codimension of the set of polynomials \ $\varphi: \C^1 \to \C^1$ \ of the form (\ref{param}) satisfying all corresponding conditions \ $\varphi(\alpha_j)=\varphi(\beta_j)$ \ or \ $\varphi'(\alpha_j)=0$ \ (if \ $\alpha_j=\beta_j$) \ in the space of all polynomials (\ref{param}). 
The {\em rank} of such a system of elementary conditions is the common value of its \ $d$-ranks for all sufficiently large \ $d$.
\end{definition}

\begin{remark} \rm
\label{rema3}
By the interpolation theorem, for this ``sufficiently large \ $d$'' \ we can take any \ $d \geq 2s-1$ \ where \ $s$ \ is the number of conditions in the system.
Accordingly, if \ $\rho < \frac{d}{2},$ \ then the rank of a system of elementary conditions is equal to \ $\rho$ \ if and only if its \ $d$-rank is equal to \ $\rho$. 
\end{remark}

\begin{definition} \rm
\label{gene}

An affine complex subspace \ $\K \subset P(d,k)$ \ is called {\em decent} if 
the closure of \ $\K$ \ in the projectivization \ $\CP^{d k} $ \ of \ $P(d,k)$ 

a) is transversal to the stratified variety consisting of the discriminant \ $\Sigma$ \ and the ``hyperplane at infinity'' \ $\CP^{d k} \setminus P(d,k)$, 

b) is transversal to the closures of all subspaces of codimension \ $k$ \ in \ $P(d,k)$ \ defined by single elementary conditions, and has non-empty intersections with each of them in \ $P(d,k)$, \ and 

c) does not intersect the closure of the set \ $\divideontimes$ \ in $\CP^{d k}$.
\end{definition}

By the generalized Lefschetz plane section theorem (see, for example, \cite{GM}, \S 2.2 of Introduction), if an affine subspace \ $\K \subset P(d,k)$ \ of complex dimension \ $D$ \ is decent, then for all \ $i < D$ \ we have
$$ H^i(P(d,k) \setminus \Sigma) \simeq H^i(\K \setminus \Sigma) \ .$$

By Alexander duality (cf. \cite{AA}), the cohomology groups of the space \ $\K \setminus \Sigma,$ \ $\dim \K =D$, \ are isomorphic to the Borel--Moore homology groups of the discriminant,
\begin{equation}
\label{alex}
\tilde H^i({\mathcal K} \setminus \Sigma) \simeq \bar H_{2D-1-i}(\Sigma \cap \K ) .
\end{equation}

It is these groups that will be calculated in the rest of the work.

\begin{proposition}
\label{compost}
1$)$ If a map \ $f = (f_1, \dots, f_k)$ \ of class \ $P(d,k)$ \ belongs to \ $\divideontimes$, \ then there exists a polynomial \ $P(z)$ \ of degree greater than $1$ and polynomials \ $Q_1, \dots, Q_k$ \ such that all \ $k$ \ components \ $f_j$ \ of \ $f$ \ have the form \ $Q_j(P(z))$. 

2$)$ If \ $f $ \ does not belong to \ $\divideontimes$ \ then it satisfies less than \ $\frac{(d-1)^2}{2}$ \ distinct elementary conditions.
\end{proposition}

\noindent
{\it Proof.} 1. If there are infinitely many points \ $x \in \C^k$ \ such that \ $f(\alpha) =x$ \ for more than one point $\alpha \in \C^1$ \ or \ $x$ \ is a critical value of \ $f$, \ then such points \ $x$ \ form a complex algebraic curve \ $C \subset \C^n$. \ Infinitely many points of \ $\C^1$ \ are mapped by \ $f$ \ to this curve, hence they all go to it. By Hartogs' theorem, this map \ $f: \C^1 \to C$ \ can be lifted to a holomorphic map \ $\check f$ \ of \ $\C^1$ \ to the normalization \ $\check C$ \ of \ $C$. 

This curve \ $\check C$ \ is simply connected. Indeed, any element of its fundamental group can be realized by a loop that avoids the critical values of \ $\check f$, \ and therefore can be lifted to a path in \ $\C^1$ \ covering our loop. Since \ $\check f$ \ is a ramified covering with finitely many preimages of each point, some finite iteration of our loop will be lifted to a closed path, which can be contracted in \ $\C^1$, \ and the projection of this contraction contracts our iterated loop in \ $\check C$. \ But the fundamental groups of complex curves have no elements of finite order greater than \ $1$, \ hence already our loop is contractible in \ $\check C$. 

Thus, \ $\check C$ \ is a simply-connected non-compact algebraic curve, hence it is isomorphic to \ $\C^1$, \ and \ $\check f$ \ is a polynomial \ $P$ \ whose degree is equal to the number of preimages of a generic point of \ $C$ \ under the map \ $f$. \ Also, the composition of the normalization map \ $\check C \to C,$ \ the identical embedding \ $C \to \C^k$, \ and the projection of \ $\C^k$ \ to the \ $j$th coordinate axis is an algebraic map \ $Q_j: \C^1 \to \C^1$ \ such that \ $f_j \equiv Q_j \circ P$, \ hence \ $Q_j$ \ is also a polynomial. 

2. At least some two components \ $f_j$, $f_{l}$ \ of such a map \ $f$ \ define a map \ $\C^1 \to \C^2$ \ satisfying only finitely many elementary conditions. 
The condition \ $f_j(\alpha) = f_j(\beta)$, $\alpha \neq \beta$, \ defines a curve of degree \ $d-1$ \ in the complex plane with coordinates \ $\alpha$ \ and \ $\beta$. \ The condition \ $(f_j -f_{l})(\alpha) = (f_j -f_{l})(\beta)$, $\alpha \neq \beta$, \ defines a curve of degree at most \ $d-2$. \ Elementary conditions satisfied by \ $f_j$ \ and \ $ f_l$ \ correspond to the intersection points of these two curves, factorized through the involution \ $(\alpha, \beta) \to (\beta, \alpha)$. \ This involution has at most \ $d-2$ \ invariant points (corresponding to the common zeros of derivatives of \ $f_j$ \ and \ $f_{l}$). \ So the number of these conditions is estimated from above by \ $\frac{(d-1)(d-2)}{2}+ \frac{d-2}{2} < \frac{(d-1)^2}{2}$. \hfill $\Box$

\begin{corollary}
\label{corcompos}
The complex dimension of the variety \ $\divideontimes \subset P(d,k)$ \ is equal to \ $a + k \frac{d}{a} - 1,$ \ where \ $a$ \ is the least divisor of \ $d$ \ greater than 1. 

The decent subspaces form a non-empty Zariski open subset in the space of all affine subspaces of dimension less than \ $d k - \dim (\divideontimes)$ \ in \ $P(d,k)$.
\hfill $\Box$
\end{corollary}

\subsection{Canonical normalization of the discriminant} 

Let \ $\K$ \ be a decent affine subspace in \ $P(d,k)$. \
For any elementary condition \ $\chi \in \X$ \ denote by \ $L(\chi)$ \ the subspace in \ $\K$ \ consisting of maps satisfying this condition. 
By item b) of Definition \ref{gene}, codimensions of all these subspaces in \ $\K$ \ are then equal to \ $k$. 

{\it Canonical normalization} \ $\widehat \Sigma$ \ of \ $\Sigma \cap \K$ \ is the subset of \ ${\X} \times {\mathcal K}$ \ consisting of all pairs \ $(\chi, f)$ \ such that \ $f \in L(\chi)$. \ The normalization map \ $\widehat \Sigma \to \Sigma \cap \K $ \ is induced by the projection \ $\X \times \K \to \K$. \ The restriction of the standard projection \ $\X \times \K \to \X$ \ to \ $\widehat \Sigma$ \ supplies \ $\widehat \Sigma$ \ with the structure of a complex affine bundle over \ $\X$ \ with fibers isomorphic to \ $\C^{\dim \K -k}$.

\subsection{Simplicial resolution of the discriminant}
\label{sire}

Let us fix a generic polynomial embedding \ $\Phi : \X \to \C^W$ \ into the space \ $\C^W$ \ of a very large dimension compared to \ $d$. \ For any finite collection of distinct points \ $\chi_1, \dots, \chi_N \in \X$, \ denote by \ $\Delta(\chi_1, \dots, \chi_N)$ \ the convex hull of all \ $N$ \ points \ $\Phi(\chi_i) \in \C^W$, \ and by \ $L(\chi_1, \dots, \chi_N)$ \ the intersection of all subspaces \ $L(\chi_1), \dots, L(\chi_N)$ \ in the decent subspace \ $\K \subset P(d,k)$. 

If the number \ $W$ \ is indeed large enough and \ $\Phi$ \ is generic, then \ $\Delta(\chi_1, \dots, \chi_N)$ \ is an \ $N$-vertex simplex for any set of points \ $\chi_i$ \ such that the space \ $L(\chi_1, \dots, \chi_N)$ \ is non-empty (by Proposition \ref{compost} (2) such numbers \ $N$ \ are uniformly bounded). Moreover, any two such simplices have only predictable intersections: these are their common faces spanning the common vertices \ $\Phi(\chi_i)$. \ We will always assume that these conditions on \ $W$ \ and \ $\Phi$ \ are satisfied.

Denote by \ $\sigma$ \ the subset of \ $\C^W \times \K$ \ equal to the union of all products 
\begin{equation}
\label{prism}
\Delta(\chi_1, \dots, \chi_N) \times L(\chi_1, \dots, \chi_N)
\end{equation}
 over all natural \ $N$ \ and all subsets \ $\{\chi_1, \dots, \chi_N\} \in B( \X,N).$ \ Denote by \ $ \Lambda$ \ the image of the projection of \ $\sigma$ \ to \ $\C^W$.

\begin{proposition}
Under the above genericity conditions on \ $\Phi$, \ the map \ $\sigma \to \Sigma \cap \K $ \ defined by the standard projection \ $\C^W \times \K \to \K $ \ is proper and surjective, and the homomorphism \ $\bar H_* (\sigma) \to \bar H_*(\Sigma) $ \ of the Borel--Moore homology groups induced by this map is an isomorphism. \hfill $\Box$
\end{proposition}

\noindent
{\it Proof} is standard, see, for example, \cite{book}, \S\S V.2.3 and III.3.4.

\subsection{Main filtration and main spectral sequence }
\label{mss}

Let \ $\K \subset P(d,k)$ \ be a decent subspace of complex dimension \ $D$. 

The resolution space \ $\sigma$ \ of \ $\Sigma \cap \K$ \ has a natural finite increasing filtration \ $\sigma_1 \subset \sigma_2 \subset \dots \subset \sigma$: \ its term \ $\sigma_\rho$ \ is the union of all products \ (\ref{prism}) \ over all systems of elementary conditions \ $\chi_1, \dots, \chi_N$ \ of \ $d$-rank \ $\rho$. \ The image \ $\Lambda \subset \C^W$ \ of the projection of \ $\sigma$ \ to \ $\C^W$ \ is also naturally filtered: its subspace \ $\Lambda_\rho$ \ is the image of \ $\sigma_\rho$.

This filtration on \ $\sigma$ \ defines a homological spectral sequence \ $E^r_{\rho,\varkappa}$ \ computing the Borel--Moore homology group of the resolution space \ $\sigma$. 

Define a cohomological spectral sequence by setting \begin{equation}
\label{aldu}
E_r^{p,q} \equiv E^r_{-p, 2D-q-1} .
\end{equation}
By Alexander duality \ (\ref{alex}), \ this spectral sequence converges to the cohomology group of \ $\K \setminus \Sigma$. 

\begin{proposition}
\label{dull}
The cohomological spectral sequence \ $($\ref{aldu}$)$ \ has non-trivial groups \ $E_1^{p,q}$ \ only in the wedge \
 $\{p < 0, q\geq -2p(k-2)\}$, \ in particular, it has only finitely many such non-zero groups on any diagonal \ $\{p+q = \mbox{const}\}$.
\end{proposition}

\noindent
{\it Proof.} 
For any \ $\rho =1, \dots, d$, \ the complex codimension of the image of the projection of \ $\sigma_\rho \setminus \sigma_{\rho-1}$ \ to \ $\K$ \ is at least \ $\rho (k - 2)$. Indeed, the space of affine planes of codimension \ $k \rho$ \ in \ $P(d,k)$, \ defined by systems of elementary conditions of \ $d$-rank \ $\rho$, \ is at most \ $2 \rho$-dimensional, so the union of the points of all these planes has (complex) codimension at least \ $\rho (k-2)$ \ in \ $P(d,k)$. By transversality condition a) of Definition \ref{gene}, the same is true for its intersection with \ $\K$. 

The preimage of any point \ $f \in \K$ \ of the image of this projection \ $\sigma_\rho \setminus \sigma_{\rho-1} \to \K$ \ is a simplex \ $\Delta(\chi_1, \dots, \chi_N) \subset \Lambda$, \ from which some of its faces (corresponding to systems of elementary conditions defining subspaces of greater dimension in \ $P(d,k)$) \ are removed. The Borel--Moore homology group of this fiber is isomorphic to the relative homology group of a certain pair of finite simplicial complexes of dimensions \ $\rho-1 $ \ and \ $\rho-2$, \ see \cite{koz} or \S VI.7.1 of \cite{book}. Namely, the first of these complexes is the order complex of all planes in \ $P(d,k)$ \ containing the point \ $f$ \ and defined by systems of elementary conditions of \ $d$-rank \ $\leq \rho$, \ and the second complex is the link of this order complex. Thus, by the Leray spectral sequence of our projection \ $\sigma_\rho \setminus \sigma_{\rho-1} \to \K$ \ the Borel--Moore homology groups of \ $\sigma_\rho \setminus \sigma_{\rho-1}$ \ are trivial in dimensions exceeding \ $2(D-\rho(k-2) )+ \rho-1$. \ By (\ref{aldu}), this implies the statement of proposition. \hfill $\Box$

\subsection{Terms \ $E^{p,q}_1$ \ for stable values of \ $p$.}

\begin{definition} \rm
\label{genegene}
A \ $D$-dimensional affine subspace \ $\K$ \ of \ $P(d,k)$ \ is {\em sufficient} if it is decent, and 
any subspace of \ $\K$ \ defined by a system of elementary conditions of rank 
\begin{equation}
\label{stabdim}
\rho < \frac{D+1}{k+2}
\end{equation}
 is non-empty and has codimension exactly \ $k \rho$ \ in \ $\K$.
\end{definition}

\begin{proposition}
\label{proest} If \ $
 D < d k - \dim (\divideontimes) $ \ 
then the sufficient subspaces form a non-empty Zariski open subset in the space of all \ $D$-dimen\-si\-onal affine subspaces of \ $P(d,k)$. 
\end{proposition}

\noindent
{\it Proof.} The set of not decent subspaces obviously is Zariski closed. By Remark \ref{rema3}, any system of elementary conditions of rank \ $\rho$ \ satisfying (\ref{stabdim}) defines a subspace of codimension \ $k \rho$ \ in the space \ $P(d,k)$. \ It is easy to calculate that 
the set of affine \ $D$-dimensional subspaces in \ $P(d,k)$ \ that are non-transversal or parallel to a particular subspace of codimension \ $k\rho$ \ has codimension \ $D-k\rho+1$ \ in the space of all affine subspaces of this dimension. The family of all subspaces of codimension \ $k \rho$ \ defined by systems of elementary conditions of rank \ $\rho$ \ is at most $2\rho$-parametric, hence the union of these sets corresponding to all such subspaces is a semialgebraic set of (complex) codimension at least \ $D+1-\rho(k+2) $ \ in the space of affine subspaces of dimension \ $D$ \ in $P(d,k)$. \ By (\ref{stabdim}), this number is positive, so the subspaces satisfying the last condition of Definition \ref{genegene} also form a non-empty Zariski open set. \hfill $\Box$ 

\begin{corollary} For any two numbers \ $d'>d$, \ if \ $\K \subset P(d,k)$ \ and \ $\K' \subset P(d',k)$ \ are sufficient subspaces of dimensions \ $D $ \ and \ $ D'$, \ then all groups \ $E_1^{p,q}(d')$, \ $p \in \left[-\min\left( \frac{D+1}{k+2}, \frac{D'+1}{k+2}  \right),-1\right]$, \ of our cohomological spectral sequence calculating the group \ $H^*(\K' \setminus \Sigma)$ \ are isomorphic to groups \ $E_1^{p,q}(d)$ \ with the same \ $p$ \ and \ $q$.
\end{corollary}

{\it Proof.} Since \ $\K$ \ is sufficient, for any \ $\rho$ \ satisfying inequality \ (\ref{stabdim}) \ the difference \ $\sigma_\rho \setminus \sigma_{\rho-1}$ \ together with the restriction of the projection \ $\C^W \times \K \to \C^W$ \ to it is a \ $(D- k \rho)$-dimensional complex affine bundle over the semialgebraic set \ $\Lambda_\rho \setminus \Lambda_{\rho-1} \subset \C^W$. \ In particular, we have the Thom isomorphism 
\begin{equation}
\label{thomthom}
E^1_{\rho, \varkappa} \equiv \bar H_{\rho + \varkappa}( \sigma_\rho \setminus \sigma_{\rho-1}) \simeq \bar H_{\rho + \varkappa- 2(D- k \rho)}(\Lambda_\rho \setminus \Lambda_{\rho-1}) \ .
\end{equation}
If \ $\K'$ \ and \ $\K$ \ are sufficient subspaces of \ $P(d',k)$ \ and \ $P(d,k)$ \ respectively, then these terms of homological spectral sequences calculating Borel--Moore homology groups of (simplicial resolutions of) spaces \ $\Sigma \cap \K'$ \ and \ $\Sigma \cap \K$ \ coincide up to the shift of index \ $q$ \ by \ $2 (D'-D)$. \ The corresponding transformations (\ref{aldu}) erase this difference. \hfill $\Box$ 
\medskip

However, to prove the commutation of these isomorphisms with higher differentials of spectral sequences we need some additional effort and additional restrictions on dimensions; this will be done in the next section.

\section{Proof of Theorem \ref{stabili}}
\label{prostab}

We assume that a number \ $k \geq 3$ \ is fixed for this entire section.

\subsection{Definition of the function \ $d_k(s)$ \ (see Theorem \ref{stabili}).}

\begin{definition} \rm
The {\em degeneracy} of a non-constant algebraic map \ $f: \C^1 \to M$, \ where \ $M$ \ is a manifold, is equal to \ $0$ \ if \ $f$ \ is a smooth embedding, otherwise it is equal to 
the sum of 
\begin{enumerate}
\item
numbers \ $a-1$ \ over all points of \ $M$ \ which are images of \ $a>1$ points of \ $\C^1$, 
\item
the number of points of \ $\C^1$ \ at which \ $f'=0$, \ and
\item the degeneracy of the induced map \ $\check f: \C^1 \to \check M$, \ where \ $\check M$ \ is the result of the blow-up of \ $M$ \ at all points of \ $M$ \ mentioned in item (1) above and all critical values of \ $f$.
\end{enumerate}
\end{definition}

The degeneracy of a map can be infinite, if the number of points mentioned in item (1) of this definition is infinite.

\begin{proposition}
\label{pro202}
For any natural \ $D$ \ there exist numbers \ $T(D)$ \ and $\delta_1(D)$ \ such that generic affine subspaces of dimensions at most \ $D$ \ in the spaces \ $P(d, k)$ \ with arbitrary \ $d \geq \delta_1(s)$ \ do not contain maps \ $f: \C^1 \to \C^k$ \ of degeneracy greater than \ $ T(D)$. \hfill $\Box$
\end{proposition}

Further, for an arbitrary map \ $f: \C^1 \to \C^k$ \ of class \ $P(d,k) \setminus \divideontimes$ \ 
consider the space \ $\mathcal J(f)$ \ of multigerms ($\psi_1, \dots, \psi_r$) \ of holomorphic maps \ $\psi_l: (\C^1, \tau_l) \to \C^k$, \ where \ $\tau_l$ \ are all critical points of \ $f$ \ and all preimages of self-intersection points of \ $f$. The restriction of \ $f$ \ to these neighborhoods is the marked point \ $\{f\}$ \ of this space. 
The group \ $\mathcal H(f)$ \ of simultaneous local holomorphic diffeomorphisms \ $\C^1 \to \C^1$, \ defined in neighborhoods of all points \ $\tau_l$, \ and local holomorphic diffeomorphisms \ $\C^k \to \C^k,$ \ defined in neighborhoods of all images of these points, acts on the space \ $\mathcal J(f)$. The notion of an {\it infinitesimally versal deformation} of this action is defined in the usual way, see \cite{AVG}. Namely, a {\it deformation} of the multigerm \ $(\psi_1, \dots, \psi_r)$ \ is a collection of maps \ $\Psi_l: (\C^1 \times \C^m) \to \C^k$ \ defined in some neighborhoods of points \ $\tau_l \times 0 \in \C^1 \times \C^m$ \ such that \ $\psi_l \equiv \Psi_l(\cdot, 0)$ \ for all \ $l$. \ So, a deformation can be considered as a family of collections of \ $r$ \ maps \ $\C^1 \to \C^k$ \ defined in neighborhoods of points \ $\tau_l$ \ and depending on the \ $m$-dimensional \ parameter \ $\mu =(\mu_1, \dots, \mu_m)$. 
Such a deformation is called {\it infinitesimally versal} if this family intersects transversally the \ $\mathcal H(f)$-orbit of the collection \ $ (\psi_1, \dots, \psi_r)$ \ in the space \ $\mathcal J(f)$ \ at the point \ $\{f\}$. \ In formal terms, this means that any element of \ $\mathcal J(f)$, \ i.e., a collection of germs \ $\theta_l(\C^1,\tau_l) \to \C^k$, $l=1, \dots, r$, \ can be represented as a sum of 

(A) a collection of maps of the form \ $V \circ \psi_l : (\C^1, \tau_l) \to \C^k$ \
where \ $V$ \ is a holomorphic vector field in \ $\C^k$ \ defined in a neighborhood of the union of all points \ $\psi_l(\tau_l)$;

(B) a collection of Lie derivatives of all maps \ $\psi_l$ \ along some holomorphic vector fields in \ $\C^1$ \ defined in neighborhoods of all points \ $\tau_l$, \ and

(C) a linear combination of the form $$\sum_{i=1}^m \alpha_i \left( \left. \frac{\partial \Psi_1}{\partial \mu_i}\right|_{\mu=0}, \dots, \left. \frac{\partial \Psi_r}{\partial \mu_i}\right|_{\mu=0} \right) , $$
where \ $\alpha_i$ \ are some complex coefficients and \ $\mu_i$ \ are parameters of the deformation \ $(\Psi_1, \dots, \Psi_r)$.

The space \ $P(d,k)$ \ provides a deformation (with parameters \ $\lambda_i$) of the collection \ $\{f\}$ \ of germs of any map \ $f \in P(d,k)$ \ at all points \ $\tau_l \in \C^1$ \ as above.

\begin{definition} \rm
A map \ $f \in P(d,k)$ \ is {\em tame} if the space \ $P(d,k)$ \ is an infinitesimally versal deformation of the action of the corresponding group \ $ \mathcal H(f)$ \ on the space \ $\mathcal J(f)$ \ at the point \ $\{f\}$.
\end{definition}

\begin{proposition}
\label{pro303}
For any natural number \ $T$ \ there exists a number \ $\delta_2(T)$ \ such that all maps \ $f \in P(d,k),$ \ $d \geq \delta_2(T)$, \ of degeneracy at most \ $T$ \ are tame. 
\end{proposition}

\noindent
{\it Proof.} Given a map \ $f \in P(d,k),$ \ let us compose it with a generic projection \ $\pi: \C^k \to \C^2$ \ and consider the set of obtained germs of plane curves at the projections of all singular points of the curve \ $f(\C^1)$. 

For any \ $T$ \ there is only a finite set of types of sets of plane curve singularities (described by Puiseux exponents and tangency orders of different local components), which can appear in this way from maps \ $f$ \ of degeneracy \ $\leq T$ \ (where \ $f$ \ can belong to spaces \ $P(d,k)$ \ with arbitrary \ $d$). 

 For each of these germs of the plane curve \ $\pi(f(\C^1))$ \ we have the following fact.

\begin{lemma} \label{lem99}
Let \ $\gamma_l: (\C^1, 0) \to (\C^2, A),$ $l=1, \dots, u$, \ be a finite collection of germs of polynomial maps such that any point of a punctured neighborhood of \ $A$ \ is the image of at most one point of the disjoint union of \ $u$ \ copies of \ $\C^1$. \ The image of the induced homomorphism from the algebra of germs of holomorphic functions \ $(\C^2, A) \to (\C^1,0)$ \ to the algebra of collections of \ $u$ \ germs of functions \ $(\C^1,0) \to (\C^1,0)$ \ $($sending any germ \ $\theta: (\C^2,A) \to (\C^1,0)$ \ to the collection of maps \ $\theta \circ \gamma_l)$ \
then contains some finite degree of the latter algebra $($i.e., there is a number \ $\nu$ \ such that any collection of \ $u$ \ functions \  $(\C^1,0) \to (\C^1,0)$  \ having zero of order \ $\geq \nu$ at the origins belongs to the image of this homomorphism$)$. This degree \ $\nu$ \ can be effectively estimated from above by the singularity type of our collection of parametric curves \ $\gamma_l$. 
\end{lemma}

\noindent
{\it Proof of  Lemma \ref{lem99}.} In the case of an irreducible curve (that is, if \ $u=1$) \ this statement easily follows from the existence of a Puiseux expansion of the corresponding map \ $\gamma_1$. \ If \ $u>1$, \ consider the similar homomorphism to the algebra of functions on the union of first \ $u-1$ \ copies of \ $\C^1$ \ (i.e. of preimages of maps \ $\gamma_1, \dots, \gamma_{u-1}$). \ By induction hypothesis this homomorphism satisfies the statement of the Lemma. Its kernel contains all holomorphic functions of the form \ $\varphi \cdot \varkappa$, \ where \ $\varphi$ \ is an arbitrary germ \ $(\C^2, A) \to (\C^1,0)$ \ and \ $\varkappa$ \ is the product of the equations of all curves \ $\gamma_l (\C^1)$, \ $l=1, \dots, u-1$. \ The restriction of \ $\varkappa$ \ to the \ $u$-th component is equal to a non-zero polynomial of its parameter; the degree of the lowest non-zero term of this polynomial is determined by the singularity type of our collection \ $\{\gamma_l\}$. \ By induction hypothesis the space of restrictions of arbitrary functions \ $\varphi$ \ contains a degree of the maximal ideal in the space of functions in this parameter, hence the same (with a greater value of the degree) is true for the space of restrictions of functions \ $\varphi \cdot \varkappa$. \hfill $\Box$

\begin{corollary} For any natural \ $T$ \ there is a number \ $\delta_3$ \ such that for any map \ $f$ \ of degeneracy \ $\leq T$ \ already the summands of type \ $($A$)$ \ from the above definition of infinitesimal versality corresponding to all singular points of the curve \ $f(\C^1)$ \ contain {\em all} collections of local maps \ $(\psi_1, \dots, \psi_r)$, \ all of whose components are defined by collections of function germs belonging to certain degrees of the maximal ideals of corresponding algebras of holomorphic germs \ $(\C^1,\tau) \to \C^1$, and the sum of these degrees over all these singular points is at most \ $\delta_3$. \hfill $\Box$ 
\end{corollary} 

By the interpolation theorem, if \ $d$ \ is large enough then the summands of type \ (C) \ for the deformation \ $P(d,k)$ \ of the collection of germs \ $\{f\}$ \ generate the quotient space of \ $\mathcal J(f)$ \ by the space of such collections \ $(\psi_1, \dots, \psi_r)$. \ This implies Proposition \ref{pro303}.
\hfill $\Box$ 

\begin{definition} \rm
\label{defidk}
Given a natural number \ $s$, \ $D(s)$ \ is the minimal natural number \ $D$ \ such that 
\begin{equation}
\label{eststa}
s \leq \min\left(D-1, \left(\left[\frac{D}{k+2}\right]+1\right)(2k-5)-2\right).
\end{equation}
The number \ $d_k(s)$ \ assumed in Theorem \ref{stabili} is equal to the maximum of the numbers \ $\delta_1(D(s))$ \ and \ $\delta_2(T(D(s)))$, \ where \ $\delta_1(\cdot)$ \ and \ $T(\cdot)$ are defined in Proposition \ref{pro202}, and \ $\delta_2(\cdot)$ \ in Proposition \ref{pro303}.
\end{definition}

\subsection{Isomorphism of spectral sequences in the stable domain}

The restriction of the main filtration of the space \ $\sigma \subset \Lambda \times \K $ \ 
to the term \ $\sigma_{\tau}$ \ of this filtration defines a spectral sequence \ $E^r_{\rho,\varkappa}(\sigma_{\tau})$ \ converging to group \ $\bar H_*(\sigma_{\tau})$. 

\begin{proposition}
\label{dull2} 
For any numbers \ $d < d'$, \ let \ $\K \subset P(d,k)$ \ and \ $\K' \subset P(d',k)$ \ be sufficient affine subspaces of dimension \ $D$, \ all points of which are tame. Then for any \ $\tau < \frac{D+1}{k+2}$ \ the 
homological spectral sequences 
defined by restricting our filtrations to the terms 
\ $\sigma_{\tau}$ \ and \ $\sigma'_{\tau}$ \ of our resolutions \ $\sigma \subset \Lambda \times \K $ \ and \ $\sigma' \subset \Lambda \times \K' $ \ respectively 
are isomorphic to each other starting from their pages \ $E_1$. 
\end{proposition}

\noindent
{\it Proof.} Both \ $\K$ \ and \ $\K'$ \ are affine subspaces of the space \ $\tilde P(d',k)$ \ of maps \ $\C^1 \to \C^k$ \ defined by systems of \ $k$ \ polynomials of the form \ $\lambda_0 z^{d'} + \lambda_1 z^{d'-1} + \dots + \lambda_{d'}.$ \ 
Our construction of simplicial resolutions cannot be applied immediately to the discriminant set \ $\Sigma \cap \tilde P(d',k)$, \ since it has points satisfying infinitely many elementary conditions; however, the analogue \ $\tilde \sigma_{\tau} \subset \Lambda_{\tau} \times \tilde P(d',k)$ \ of the \ $\tau$-th term of this simplicial resolution can be constructed in exactly the same way as previously.

Since all points \ $f$ \ of \ $\K$ \ are tame, the subspace \ $P(d,k)$ \ of \ $\tilde P(d',k)$ \ is transversal to the stratified variety \ $\Sigma \cap \tilde P(d',k)$ \ at all points of its subspace \ $\K$. \ Since \ $\K$ \ is a decent subspace of \ $P(d,k)$ \ (see Definition \ref{gene}), \ it also is transversal to \ $\Sigma \cap \tilde P(d',k)$ \ in \ $\tilde P(d',k)$. \
Let \ {\large $B$} \ be a huge ``exhausting'' open ball in \ $\tilde P(d',k)$ \ such that the identical embeddings induce isomorphisms \ $\bar H_*(X\cap \mbox{\large $B$}) \simeq \bar H_*(X)$ \ for all involved algebraic varieties \ $X$ \ (such as  \ $\Sigma \in \tilde P(d',k)$, $\Sigma \cap P(d,k)$, $\Sigma \cap \K$, \ different their strata, etc), and the boundary of \ {\large $B$} \ is transversal to all these varieties.

 By Thom isotopy lemma (see e.g. \cite{GM}, \S I.1.5), there exists a tubular neighborhood \ $U$ \ of the subspace \ $\K$ \ in \ $\tilde P(d,k) \cap \mbox{\large B}$ \ such that the pair \ $(U, \Sigma \cap U)$ \ is homeomorphic to the direct product of the pair \ $(\K \cap \mbox{\large $B$}, \Sigma \cap \K \cap \mbox{\large $B$})$ \ and an open ball of dimension \
$2 (\dim_{\C}\tilde P(d',k)-\dim_{\C} \K) \equiv 2((d'+1)k-D)$. \ Consider the subset \ $\tilde \sigma_{\tau}(U)$ \ of the space \ $\tilde \sigma_{\tau} \subset \Lambda_\tau \times \tilde P(d',k)$ \ consisting of only the points whose projections to \ $\tilde P(d',k)$ \ belong to \ $U$. \ It can be considered as the \ $\tau$-th term of the simplicial resolution of the set \ $\tilde \Sigma(d',k) \cap U$. \ By construction, it is also homeomorphic to the product of \ $\sigma_{\tau}$ \ and an open ball of dimension \ $2 ((d'+1)k-D)$. \
In particular, the spectral sequences calculating the Borel--Moore homology groups of these spaces are isomorphic to each other up to a shift of dimensions: \
$E^r_{\rho, \varkappa}(\tilde \sigma_\tau(U)) \simeq E^r_{\rho,\varkappa-2((d'+1)k-D)}(\sigma_\tau)$ \
for all \ $r \geq 1$, \ $\rho \leq \tau$ \ and any \ $\varkappa$. \ On the other hand, the identical embedding \ $\tilde \sigma_\tau(U) \to \tilde \sigma_\tau$ \ induces a homomorphism of the corresponding spectral sequences. This homomorphism is an isomorphism of all terms \ $E^1_{\rho, \varkappa}$: \ indeed, any set \ $\tilde \sigma_\rho \setminus \tilde \sigma_{\rho-1}$ \ is the space of an affine bundle with base \ $\Lambda_\rho \setminus \Lambda_{\rho-1}$, \ and its subset \ $\tilde \sigma_\rho(U) \setminus \tilde \sigma_{\rho-1}(U)$ \ is the space of a fiber bundle, the base of which is a subspace of \ $\Lambda_\rho \setminus \Lambda_{\rho-1}$ \ having the same Borel--Moore homology groups, and the fibers are open balls in the fibers of the former bundle. Thus, our spectral sequences calculating the Borel--Moore homology groups of spaces \ $\sigma_\tau \subset \Lambda_\tau \times \K$ \ and \ $\tilde \sigma_\tau \subset \Lambda_\tau \times \tilde P(d',k)$ \ are isomorphic up to a shift of indices \ $\varkappa$ \ by \ $2((d'+1)k-D)$.

In the same way we prove that our spectral sequence computing the homology groups of \ $\sigma'_\tau$ \ is isomorphic to the same spectral sequence for \ $\tilde \sigma_\tau$ \ up to the same shift of indices \ $\varkappa$. \ In particular, our spectral sequences for the spaces \ $\sigma_\tau$ \ and \ $\sigma'_\tau$ \ are isomorphic to each other.
\hfill $\Box$ \medskip

\begin{corollary}
\label{cordul}
1. In conditions of Proposition \ref{dull2}, spectral sequences \ $($\ref{aldu}$)$ \ converging to cohomology groups of corresponding spaces \ $\K' \setminus \Sigma'$ \ and \ $\K \setminus \Sigma$ \ of complex knots are isomorphic to each other starting from term \ $E_1$ \ in the domain of the \ $(p, q)$-plane where 
\begin{equation}
\label{e12}p+q \leq \left(\left[\frac{D}{k+2}\right]+1\right)(2k-5)-2.
\end{equation} 
In particular, groups \ $H^i(\K \setminus \Sigma)$ \ and \ $H^i(\K' \setminus \Sigma)$ \ are isomorphic to each other for \ $i$ \ not exceeding the right-hand part of \ $($\ref{e12}$)$.

2. The last isomorphisms are natural.
\end{corollary}

\noindent
{\it Proof.} 1. By Proposition \ref{dull} all non-zero terms \ $E_1^{p,q}$ \ of these cohomological spectral sequences with \ $p < -\left[\frac{D}{k+2}\right]$ \ (which only can be different for these two spectral sequences or provide different differentials) lie in the domain of the \ $(p,q)$-plane, where \ $p+q \geq \left(\left[\frac{D}{k+2}\right]+1\right)(2k-5) $. 

2. We can connect our subspaces \ $\K$ \ and \ $\K'$ \ by a path \ $\{\K_\tau\}$, $\tau \in [0,1]$, \ in the space of affine \ $D$-dimensional subspaces of \ $\tilde P(d',k)$ \ satisfying all the same genericity conditions. Consider the space of pairs \ $(\tau, f)$ \ where \ $\tau \in [0,1]$ \ and $f \in \K_\tau \setminus \Sigma$. The inclusion of any fiber \ $\K_\tau \setminus \Sigma$ \ to this space induces then an isomorphism of all homology groups in dimensions not exceeding the right-hand part of (\ref{e12}). \hfill $\Box$

\begin{corollary}
\label{cor15}
Under the conditions of Proposition 13, the groups \ $H^i(P(d,k) \setminus \Sigma)$ \ and \ $H^i(P(d',k) \setminus \Sigma)$ \ are naturally isomorphic to each other for \ $i$ \ not exceeding the right-hand part of \ $($\ref{eststa}$)$.
\end{corollary}

\noindent
{\it Proof} of this corollary is the composition of Corollary \ref{cordul} and strong Lefschetz plane section theorem, see \cite{GM}. \hfill $\Box$ \medskip

\noindent
{\it Proof of Theorem \ref{stabili}}. By Propositions \ref{pro202} and \ref{pro303}, for any natural \ $s$ \ there are numbers \ $d_k(s)$ \ and \ $D(s)$ \ (see Definition \ref{defidk}) such that the conditions of Proposition \ref{dull2} are satisfied for generic \ $D$-dimensional subspaces \ of spaces \ $P(d,k)$ \ with \ $d\geq d_k(s)$ \ and \ $D=D(s)$. Then Theorem \ref{stabili} follows from Corollary \ref{cor15}. \hfill $\Box$

\begin{table}
\caption{Column $p=-3$ of main spectral sequence}
\label{e3}
\begin{tabular}{|c||c|c|c|c|c|c|c|}
\hline
$q$ & $6k-11$ & $6k-10$ & $6k-9$ & $6k-8$ & $6k-7$ & $6k-6$ & $6k-5$ \\
\hline
$E_1^{-3,q}$ & $\Z_2$ & $\Z_3$ & $\Z \oplus \Z_3$ & $ T$ & $ T$ & $\Z \oplus T$ & $T$ \\
\hline
\end{tabular}

\end{table}

\unitlength 0.7mm

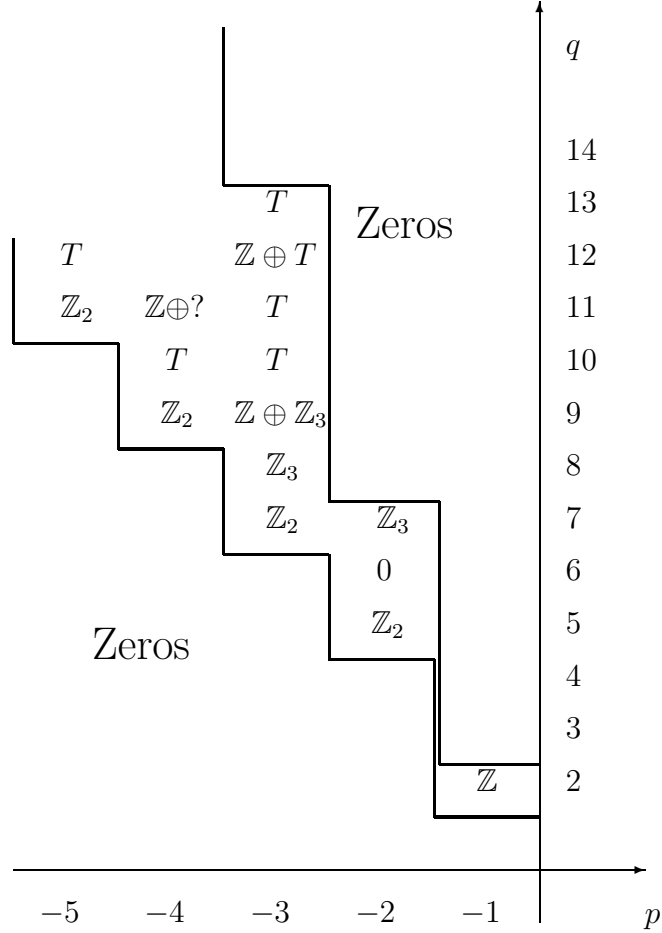
\begin{figure}

\begin{picture}(120,175)
\put(105,165){$q$}
\put(0,10){\vector(1,0){120}}
\put(5,0){$-5$}
\put(25,0){$-4$}
\put(45,0){$-3$}
\put(65,0){$-2$}
\put(85,0){$-1$}
\put(100,0){\vector(0,1){175}}
\put(120,0){$p$}
\put(105,25){$2$}
\put(105,35){$3$}
\put(105,45){$4$}
\put(105,55){$5$}
\put(105,65){$6$}
\put(105,75){$7$}
\put(105,85){$8$}
\put(105,95){$9$}
\put(105,105){$10$}
\put(105,115){$11$}
\put(105,125){$12$}
\put(105,135){$13$}
\put(105,145){$14$}
\put(88,25){$\Z$}
\put(68,55){$\Z_2$}
\put(69,65){$0$}
\put(69,75){$\Z_3$}
\put(48, 75){$\Z_2$}
\put(48,85){$\Z_3$}
\put(42,95){\bf $\Z \oplus \Z_3$}
\put(48,105){$T$}
\put(48,115){$T$}
\put(42,125){$\Z \oplus T$}
\put(48,135){$T$}
\put(28,95){$\Z_2$}
\put(29,105){$T$}
\put(25,115){$\Z \oplus ?$}
\put(9,115){$\Z_2$}
\put(9,125){$T$}
\put(15,50){\Large{Zeros}}
\put(65,130){\Large{Zeros}}
\thicklines
\put(0,110){\line(0,1){20}}
\put(0,110){\line(1,0){20}}
\put(20,110){\line(0,-1){20}}
\put(20,90){\line(1,0){20}}
\put(40,90){\line(0,-1){20}}
\put(40,70){\line(1,0){20}}
\put(60,70){\line(0,-1){20}}
\put(60,50){\line(1,0){20}}
\put(80,50){\line(0,-1){30}}
\put(80,20){\line(1,0){20}}
\put(40,140){\line(0,1){30}}
\put(40,140){\line(1,0){20}}
\put(60,140){\line(0,-1){60}}
\put(60,80){\line(1,0){21}}
\put(81,80){\line(0,-1){50}}
\put(81,30){\line(1,0){19}}
\end{picture}
\caption{Page \ $E_1$ \ of main spectral sequence for \ $k=3$}
\label{sps0}
\end{figure}

\subsection{Main technical result}

Also by Propositions \ref{pro202} and \ref{pro303}, and Corollary \ref{cordul}, when \ $d$ \ grows to infinity, our spectral sequences stabilize to a universal (depending only on \ $k$) \ spectral sequence. 
 We will call it the {\em main cohomological spectral sequence} of our problem. 

\begin{theorem}[see Fig.~\ref{sps0}]
\label{mthm2}
1. The column \ $E_1^{-1,*}$ \ of main cohomological spectral sequence contains only one non-zero group \ $E_1^{-1,2k-4} \sim \Z$. 

2. The column \ $E_1^{-2,*}$ \ of this spectral sequence contains only two non-zero groups, \ $E_1^{-2,4k-7} \sim \Z_2$ \ and \ $E_1^{-2,4k-5} \sim \Z_3$. 

3. All groups \ $E_1^{-3,q}$ \ with \ $q \not \in [6k-11,6k-5]$ \ are trivial, while such groups with \ $q \in [6k-11,6k-5]$ \ are as shown in Table \ref{e3}, where \ $T$ \ denotes finite $($in general, different and possibly trivial$)$ groups. 

4. Group \ $E_1^{-4,8k-13}$ \ is infinite.

5. For any \ $p \leq -2$, \ group \ $E_1^{p,q}$ \ of main spectral sequence 

\ \ a$)$ is trivial for \ $q \leq -p(2k-4) $, 

\ \ b$)$ is isomorphic to \ $\Z_2$ \ for \ $q=-p(2k-4)+1$, \ and 

\ \ c$)$ is finite for \ $q \leq -p(2k-4)+2$.
\end{theorem}

The picture of this spectral sequence for \ $k > 3$ \ can be obtained from Fig.~\ref{sps0} by the shift of any column \ $\{E_1^{p,*}\}$ \ by \ $-2p(k-3)$ \ in the vertical direction.

Statements 1, 2, 3, 4 and 5 of this theorem will be proved in Sections \ref{frst}, \ref{scnd}, \ref{sec3t}, \ref{prof4} and \ref{prof5}, respectively. Theorem \ref{mthm} follows directly from this one.

\section{On the filtration terms in the stable range}

\label{prelim}

Spaces \ $\Lambda_\rho \setminus \Lambda_{\rho-1}$ \ for \ $\rho$ \ in the stable range (i.e. satisfying (\ref{stabdim})) are naturally stratified according to the structure of systems of elementary conditions of rank \ $\rho$.

\begin{example} \rm
\label{ex1}
Some such systems of rank 2 are pairs of conditions \ $f(\alpha)=f(\beta), f(\gamma)=f(\delta),$ \ where all points \ $\alpha, \beta, \gamma, \delta$ \ are distinct; some others are of the form \ $f(\alpha)=f(\beta)=f(\gamma)$. The corresponding simplex \ $\Delta(\cdot)$ \ in \ $\Lambda_2$ \ in the first case is a segment (whose endpoints \ $\Phi(\{\alpha, \beta\})$ \ and \ $\Phi(\{\gamma, \delta\})$ \ lie in \ $\Lambda_1$), while in the second case we have four such simplices in \ $\Lambda_2 \setminus \Lambda_1$: \ the triangle spanned by points \ $\Phi(\{\alpha, \beta\})$, $\Phi(\{\beta, \gamma\})$, $\Phi(\{\gamma, \alpha\})$, \ and each of its three edges (whose endpoints again belong to $\Lambda_1$).
\end{example}

\begin{definition}[see \cite{ks}] \rm 
Let \ $A = (a_1, \dots, a_s)$ \ be an unordered set of natural numbers (some of which may be the same), all of whose elements \ $a_j$ \ are greater than 1. Then an \ $A$-{\it configuration} is an arbitrary set of \ $a_1 + \dots + a_s$ \ pairwise distinct points of \ $\C^1$ \ divided into subsets of cardinalities \ $a_1, \dots, a_s$. \ If additionally \ $b $ \ is a non-negative integer, then a $(A,b)$-{\it configuration} in \ $\C^1$ \ is an arbitrary \ $A$-configuration complemented by \ $b$ \ pairwise distinct points (some of which may coincide with points of the \ $A$-configuration). 
Number
\begin{equation}
(a_1-1)+ \dots + (a_s-1) +b
\label{comple}
\end{equation}
is called the {\it complexity} of any \ $(A,b)$-configuration. 
\end{definition}

If a \ $D$-dimensional subspace \ $\K \subset P(d,k)$ \ is sufficient, then for any number \ $\rho$ \ satisfying (\ref{stabdim}) there is a one-to-one correspondence between \ $(A,b)$-configurations of complexity \ $\rho$ \ and subspaces of codimension \ $k \rho$ \ in \ $\K$ \ defined by systems of elementary conditions of rank \ $\rho$. \ Namely, any \ $(A,b)$-configuration defines a subspace consisting of maps that glue points of any of the \ $s$ \ subsets of its \ $A$-configuration and have zero derivative at any of additional \ $b$ \ points.

\subsection{Simplices associated with $(A,b)$-configurations and complexes of connected graphs} 
\label{congr}
Any \ $(A,b)$-configuration \ $\Gamma,$ \ $A=(a_1, \dots, a_s)$, \ also defines a simplex \ $\Delta(\Gamma)$ \ with \ $\binom{a_1}{2} + \dots + \binom{a_s}{2} + b$ \ vertices in \ $\C^W$: \ this spans all points \ $\Phi(\{\alpha,\beta\})$, \ where \ $\alpha$ \ and \ $\beta$ \ are some two points of one of \ $s$ \ subsets of this configuration, and all points \ $\Phi(\{\alpha, \alpha\})$ \ such that \ $\alpha$ \ is a point of its \ $b$-part. Thus, the term \ $\Lambda_\rho $ \ of main filtration of \ $\Lambda$ \ is the union of such simplices defined by all \ $(A,b)$-configurations of complexity \ $ \leq \rho$.

If the symbol \ $A$ \ consists of a single number, \ $A = (a_1)$, \ then the faces of such a simplex associated with any \ $A$-configuration are in a natural one-to-one correspondence with simple graphs with \ $a_1$ \ vertices corresponding to points of this configuration: we draw an edge connecting some two vertices \ $\alpha$ \ and \ $\beta$ \ if and only if the point \ $\Phi(\{\alpha, \beta\})$ \ is one of the vertices of this face. This simplex belongs to the term \ $\Lambda_{a_1-1}$ \ of the main filtration, moreover, all its faces corresponding to not connected graphs (in particular, graphs with isolated vertices) belong to \ $\Lambda_{a_1-2}$. \ Similarly, the faces of a simplex associated with an arbitrary \ $(A,b)$-configuration are characterized by collections of \ $s$ \ graphs on \ $a_1, \dots , a_s$ \ vertices, and by additional marking or not each of the \ $b$ \ singular points. This entire simplex lies in \ $\Lambda_\rho$ \ where \ $\rho$ \ is the complexity \ $\sum (a_i-1)+b$ \ of the symbol \ $(A,b)$. \ Interior points of a face do not belong to the lower term \ $\Lambda_{\rho-1}$ \ of filtration if and only if all corresponding \ $s$ \ graphs are connected and all \ $b$ \ singular points are marked.

Thus, the Borel--Moore homology groups of the parts of these simplices lying in \ $\Lambda_\rho \setminus \Lambda_{\rho-1}$ \ are described in the following terms.

\begin{definition}
\label{congrdef}
The {\em complex of connected graphs} on \ $a$ \ vertices is the factor complex of the simplicial complex generated by faces of the \ $\binom{a}{2}$-vertex simplex by the subcomplex generated by faces corresponding to non-connected graphs.
\end{definition}

\begin{proposition}[see e.g. \cite{book}, \S V.3]
\label{congrprop}
The complex of connected graphs on \ $a$ \ vertices is acyclic in all dimensions other than \ $a-2$, \ and its \ $(a-2)$-dimensional homology group is isomorphic to \ $\Z^{(a-1)!}$.
\end{proposition}

\begin{corollary}
For any \ $(A,b)$-configuration \ $\Gamma$ \ of complexity \ $\rho$, \ $A=(a_1, \dots, a_s)$, \ the Borel--Moore homology group \ $\bar H_i(\Delta(\Gamma) \setminus \Lambda_{\rho-1})$ \ of the part of \ $\Delta(\Gamma)$ \ lying in \ $\Lambda_\rho \setminus \Lambda_{\rho-1}$ \ is trivial in all dimensions other than \ $\rho-1$ \ and is free Abelian of rank \ $\prod_{m=1}^s (a_m-1)!$ \ for \ $i=\rho-1$. 
\end{corollary}

\noindent
Indeed, the whole simplex \ $\Delta(\Gamma)$ \ can be considered as the join of its \ $ s$ \ faces spanned by the vertices \ $\Phi(\alpha, \beta)$ \ where \ $\alpha$ \ and \ $\beta$ \ belong to the same subset of the configuration, and additional \ $b$ \ points \ $\Phi(\alpha, \alpha)$ \ where \ $\alpha$ \ belongs to the \ $b$-part of the configuration. The set \ $\Delta(\Gamma) \cap (\Lambda_\rho \setminus \Lambda_{\rho-1})$ \ is the
union of interior parts of faces of \ $\Delta(\Gamma)$ \ which are joins of the parts of these faces corresponding to connected graphs, and of all \ $b$ \ additional points. \hfill $\Box$

\subsection{Inverse auxiliary filtration (see \cite{twocon})} \label{aux}

 There is a convenient filtration \ $\Theta_0 \subset \Theta_1 \subset \dots \subset \Theta_{\rho-1} = \Lambda_\rho \setminus \Lambda_{\rho-1}$ \ in any space \ $\Lambda_\rho \setminus \Lambda_{\rho-1}$. \ For example, if \ $\rho = 2$ \ then the subspace \ $\Theta_0$ \ consists of all intervals from the first case considered in Example \ref{ex1} and the edges of triangles from the second case; the set \ $\Theta_1 \setminus \Theta_0$ \ consists only of interior parts of these triangles. For certain historical reasons, this filtration is called the {\em inverse auxiliary filtration} or just {\em inverse filtration} of \ $\Lambda_\rho \setminus \Lambda_{\rho-1}$.

\begin{definition} \rm 
\label{aux1}
For any symbol \ $A= (a_1 , a_2 , \dots , a_s)$ \ as above, denote by \ $|A|$ \ the sum \ $a_1 + \dots + a_s$ \ and by \ $\#(A)$ \ the number \ $s$ \ of elements \ $a_j$ \ in \ $A$.

 The {\it defect} of an \ $A$-configuration is equal to twice its complexity minus the number of geometrically distinct points of this configuration (obviously it is a non-negative integer, and can be defined also as the difference of the complexity and \ $\#(A)$). \
The term \ $\Theta_j$ \ of the auxiliary inverse filtration of the space \ $\Lambda_\rho \setminus \Lambda_{\rho-1} $ \ is defined as the closure in this space of the union of simplices \ $\Delta(\Gamma)$ \ defined by \ $A$-configurations \ $\Gamma$ \ of complexity \ $\rho$ \ and defect \ $\leq j$. 
\end{definition}

Let us reveal the operation of closure in this definition. 

First, any simplex \ $\Delta(\Gamma)$ \ in \ $\Lambda$ \ defined by a \ $(A,b)$-configuration with \ $b >0$ \ belongs to the closure of the set of similar simplices defined by \ $A'$-configurations where symbol \ $A'$ \ is obtained from \ $A$ \ by adding \ $b$ \ numbers \ $2$. \ Therefore the entire of \ $\Lambda_\rho \setminus \Lambda_{\rho-1}$ \ is indeed covered by the closures of terms \ $\Theta_j$.

Further, for any symbol \ $A=(a_1, \dots, a_{\#(A)})$ \ define the corresponding configuration space \ $\B(A)$ \ as the 
space of unordered collections of \ $|A|$ \ points in \ $\C^1$ \ (some of which may be the same) split into subcollections of cardinalities \ $a_1, \dots, a_{\#(A)}$ \ such that all points of any subcollection of cardinality \ $a_j >2$ \ are pairwise distinct. 

To fix a topology on the space \ $\B(A)$, \ we realize it as the Cartesian product of the spaces \ $B(\C^1,a_j)$ \ over all indices \ $j=1, \dots, \#(A)$ \ with \ $a_j>2$, \ and spaces \ $\X$ \ corresponding to all \ $a_j=2$, \ factorized through permutations of such factors corresponding to equal values of \ $a_j$. \ The points of \ $\B(A)$ \ are called \ $\bar A$-configurations.

Any \ $\bar A$-configuration \ $\Gamma$ \ defines an affine subspace \ $L(\Gamma)$ \ of the space of polynomials (\ref{param}): it consists of maps that take equal values at all points of each of \ $\#(A)$ subcollections, and have zero derivative at all points \ $\alpha \in \C^1$ \ such that our \ $\bar A$-configuration \ $\Gamma$ \ contains the subcollection \ $\{\alpha, \alpha\}$. \ The codimension of this subspace \ $L(\Gamma)$ \ does not depend on \ $d$ \ if \ $d$ \ is sufficiently large (say, \ $d > |A|$). \ This stable codimension does not exceed \ $|A|-\#(A)$; denote by \ $\tilde \B(A)$ \ the Zariski open subset of \ $\B(A)$ \ consisting of configurations for which this codimension is exactly equal to \ $|A|-\#(A)$.

\begin{lemma}
\label{lem100} 
An \ $\bar A$-configuration belongs to the space \ $\tilde \B(A)$ \ if and only if 

1$)$ no two of its subcollections have a pair of common distinct points, 

2$)$ it has no two equal subcollections of the form \ $(\alpha, \alpha)$, \ and

3$)$ there are no closed chains of \ $\geq 3$ \ its subcollections such that any two neighboring subcollections of the chain have a common point. 
\end{lemma}

\noindent This lemma is elementary. \hfill $\Box$
\medskip

Like \ $A$-configurations, any \ $\bar A$-configuration \ $\Gamma$ \ defines a polyhedron \ $\Delta(\Gamma)$ \ in \ $\C^W$, \ namely, it is the convex hull of all points \ $\Phi(\{\alpha, \beta\})$, \ where \ $\alpha$ \ and \ $\beta$ \ are two distinct points of one of subcollections of \ $\Gamma$, \ and all points \ $\Phi(\{\alpha,\alpha\})$ \ over subcollections of the form \ $\{\alpha,\alpha\}$. \ If \ $\bar A$-configuration \ $\Gamma$ \ belongs to \ $\tilde \B(A)$ \ (and the map \ $\Phi$ \ satisfies genericity conditions of \S \ref{sire}), then \ $\Delta(\Gamma)$ \ is a simplex with \ $\binom{a_1}{2} + \dots + \binom{a_{\#(A)}}{2}$ \ vertices. 
Now, a new definition of main and inverse filtrations on \ $\Lambda$ \ can be formulated.

\begin{definition} 
\label{aux2}
\rm 
$\bar A$-configurations satisfying three restrictions of Lemma \ref{lem100} are called {\em regular}. 
For any symbol \ $A=(a_1, \dots, a_{\#(A)})$ \ as above, the \ {\it $A$-block} in \ $\Lambda \subset \C^W$ \ is the union of simplices \ $\Delta(\Gamma)$ \ over all regular \ $\bar A$-configurations \ $\Gamma$. \ For any natural \ $\rho$, \ the \ $\rho$-th term \ $\Lambda_\rho$ \ of the {\em main filtration} of the space \ $\Lambda$ \ is the union of \ $A$-blocks over all symbols \ $A$ \ with \ $|A|-\#(A) \leq \rho$. 
For any natural \ $\rho$ \ and any \ $j \in \{0, 1, \dots, \rho-1\}$, the \ $j$-th term \ $\Theta_j$ \ of the {\em inverse filtration } of the term \ $\Lambda_\rho \setminus \Lambda_{\rho-1}$ \ of the main filtration of \ $\Lambda$ \ is the intersection of this term \ $ \Lambda_\rho \setminus \Lambda_{\rho-1}$ \ and 
the union of \ $A$-blocks over all symbols \ $A$ \ with \ $|A|-\#(A)=\rho$ \ and $|A|-2\#(A) \leq j$.
\end{definition}

\begin{lemma}
Definition \ref{aux2} of the main and inverse filtrations is equivalent to their definitions given in \S \ref{mss} and Definition \ref{aux1}, respectively.
\end{lemma}

\noindent
{\it Proof} \ is straightforward: the expansion of the notion of \ $A$-con\-fi\-gu\-rations to that of \ $\bar A$-configurations is just an implementation of the word ``closure'' in Definition \ref{aux}, cf. \cite{twocon}. \hfill $\Box$ \medskip

As in \S \ref{congr}, any face of the simplex \ $\Delta(\Gamma)$ \ defined by a \ $\bar A$-configuration \ $\Gamma$ \ is characterized by a system of \ $\#(A)$ \ graphs on \ $a_1, \dots, a_{\#(A)}$ \ vertices associated with points of subcollections of \ $\Gamma$; \ in particular if some subcollection of \ $\Gamma$ \ is of type \ $\{\alpha, \alpha\}$, \ then the corresponding graph is a disjoint pair of points
 (respectively, two points together with a segment connecting these points) if the point \ $\Phi(\{\alpha, \alpha\})$ \ is not (respectively, is) a vertex of the face. 

\begin{definition} \rm
\label{twocc}
A simple graph is {\em two-connected} if it is connected, and removing from it an arbitrary its vertex with all incident edges, we again obtain a connected graph. The complex of two-connected graphs on a given set of vertices is defined analogously to Definition \ref{congrdef}, only with replacement of faces corresponding to connected graphs by those corresponding to two-connected graphs.
\end{definition}

\begin{theorem}[see \cite{BBB}, \cite{tur}]
\label{2c}
The complex of two-connected graphs on \ $a$ \ vertices is acyclic in all dimensions other than \ $2a-4$, \ and its \ $(2a-4)$-dimensional homology group is isomorphic to \ $\Z^{(a-2)!}$. \hfill $\Box$
\end{theorem}

\begin{lemma}[see \cite{twocon}]
\label{lem18}
For any natural \ $\rho$ \ and any \ $j \in \{0, 1, \dots, \rho-1\}$, 
term \ $\Theta_j \setminus \Theta_{j-1} $ \ of inverse filtration of the term \ $\Lambda_\rho \setminus \Lambda_{\rho-1}$ \ of main filtration of \ $\Lambda$ \ is the disjoint union of certain subsets of all \ $A$-blocks such that \ $|A|-\#(A)=\rho$ \ and \ $|A|- 2\#(A) = j$. \ Namely, any such subset of any such \ $A$-block is a fiber bundle, the base of which is the corresponding configuration space \ $\tilde \B(A)$, \ and the fiber over any regular \ $\bar A$-configuration \ $\Gamma \in \tilde \B(A)$ \ is equal to the
 union of interior points of all faces of the simplex \ $\Delta(\Gamma)$, \ such that all corresponding \ $\#(A)$ \ graphs are two-connected. \hfill $\Box$
\end{lemma}

\begin{corollary}
\label{cor1}
For any regular \ $\bar A$-configuration \ $\Gamma$, 
the Borel--Moore homology group of the fiber over the point \ $\{\Gamma\}$, \ described in the end of Lemma \ref{lem18}, is nontrivial only in dimension \ 
$2|A|-3\#(A)-1$ \ and is isomorphic to a free Abelian group of rank \ $\prod_{j=1}^{\#(A)} (a_j-2)! $. \hfill $\Box$
\end{corollary}

\begin{remark} \rm
In further calculations, we use the same notation \ $\Theta_0, \Theta_1$, etc. for terms of inverse filtrations of spaces \ $\Lambda_\rho \setminus \Lambda_{\rho-1}$ \ with different \ $\rho$: I hope this will not cause any confusion.
\end{remark}

\section{First term of the main filtration}
\label{frst}

This term \ $\sigma_1$ \ is isomorphic to the canonical normalization \ $\widehat \Sigma$ \ of the discriminant, i.e. the space of a fiber bundle over \ $\Lambda_1 \equiv \Phi(\X) \simeq \C^2$ \ with fibers equal to complex affine subspaces of codimension \ $k$ \ in \ $ {\mathcal K}$. Thus group \ $E^1_{1,\varkappa} $ \ of the main homological spectral sequence is isomorphic to \ $\Z$ \ for \ $\varkappa +1=2D-2k+4$ \ and is trivial in all other dimensions. By virtue of (\ref{aldu}), the corresponding column \ $p=-1$ \ of the main {\em cohomological} spectral sequence is as described in the first statement of Theorem \ref{mthm2}.

\section{The second term}
\label{scnd}

Statement 2) of Theorem \ref{mthm2} concerning column \ $E_1^{-2,*}$ \ follows directly from identity (\ref{aldu}), the Thom isomorphism (\ref{thomthom}) for \ $\rho=2$, \ and the following calculation of the Borel--Moore homology group \ $\bar H_*(\Lambda_2 \setminus \Lambda_1).$

\begin{theorem}
\label{thmf2}
 Let \ $\{\E^r_{j, \mathfrak{q}} \}$ \ be the spectral sequence calculating group 
\ $\bar H_*(\Lambda_2 \setminus \Lambda_1)$ \ 
and generated by the auxiliary inverse filtration \ $\{\Theta_0 \subset \Theta_1\}$ \ of \ $\Lambda_2 \setminus \Lambda_1$. 
Then all non-trivial groups of its term \ $\E^1$ are shown in Table \ref{tab2} $($left$)$, and its differential \ $\partial^1: \E^1_{1,6} \to \E^1_{0,6}$ \ is an isomorphism.
\end{theorem}

\begin{table}
\caption{Page \ $\E^1$ \ of inverse spectral sequences for $\Lambda_2 \setminus \Lambda_1$ \ and \ $\Lambda_3 \setminus \Lambda_2$}
\label{tab2}
\begin{tabular}{c|c|c|r}
$\mathfrak{q}$ & & &\\
& & &\\
8 & $\Z_2$ & 0 & \\
7 & 0 & 0 & \\
6 & $\Z_2$ & $\Z_2$ & \\
5 & 0 & $\Z_3$ & \\
4 & 0 & 0 & \\
\hline
& 0 & 1 & $j$
\end{tabular} \hspace{3cm}
\begin{tabular}{c|c|c|c|r}
$\mathfrak{q}$ & & & & \\
& & & & \\
13 & $\Z_2$ & 0 & 0 & \\
12 & $\Z_3$ & 0 & 0 & \\
11 & $\Z_2$ & $\Z_2$ & 0 & \\
10 & 0 & $\Z_3$ & 0 & \\
9 & 0 & 0 & $\Z$ & \\
8 & $\Z_6$ & $\Z \oplus \Z_2$ & $\Z \oplus \Z_2$ & \\
7 & $\Z_3$ & $\Z$ & $\Z_2$ & \\
\hline
& 0 & 1 & 2 & $j$
\end{tabular}

\end{table}

\noindent
{\it Proof.} The term \ $\Theta_0$ \ in this case is a fiber bundle with base \ $B(\X,2)$; \ its fiber over a pair of points \ $\{\chi_1 \neq \chi_2\} \subset \X$ \ is the interval in \ $\C^W$ \ with endpoints \ $\Phi(\chi_1)$ \ and \ $\Phi(\chi_2)$ \ (these endpoints lie in \ $\Lambda_0$). \ So we have 
\begin{equation}
\bar H_i(\Theta_0) \simeq \bar H_{i-1}(B(\X,2), \pm \Z)
\label{thomfour}
\end{equation}
 for any \ $i$. 
Since \ $\X \simeq \R^4$, \ the space \ $B(\X,2)$ \ is homeomorphic to \ $\R^4 \times (0,\infty) \times \RP^3$, \ and the group \ $\bar H_*(B(\X,2), \pm \Z)$ \ is isomorphic to \ $\Z_2$ \ in dimensions $5$ and $7$ and is trivial in all other dimensions. This gives us the column \ $j=0$ \ of the spectral sequence shown in Table \ref{tab2} (left). 

Remaining part \
$\Theta_1 \setminus \Theta_0$ \ of \ $\Lambda_2 \setminus \Lambda_1$ \
is a fiber bundle over \ $B(\C^1,3)$: \ its fibers are open triangles, the orientation of which is changed by the monodromy over loops in the base that define odd permutations of three points. So \ $\bar H_N(\Theta_1 \setminus \Theta_0) \simeq \bar H_{N-2}(B(\C^1,3),\pm \Z)$. 

\begin{lemma}
\label{lem6}
Group \ $\bar H_i(B(\C^1,3), \pm \Z)$ \ is isomorphic to \ $\Z_2$ \ for \ $i= 5,$ \ to \ $\Z_3$ \ for \ $i=4,$ \ and is trivial for all other \ $i$. 
\end{lemma}

\noindent
This easily follows from calculations in the standard (introduced in \cite{fuks}) cell decomposition of \ $B(\C^1,3)$, \ see e.g. \cite{book}, \S I.4. \hfill $\Box$ \medskip

This lemma gives us column \ $j=1$ \ of Table \ref{tab2} (left). Let us calculate its differential \ $\partial^1: \E^1_{1,6} \to \E^1_{0,6}$.

Group \ $\E^1_{1,6}$ \ of this spectral sequence is generated by the homology class \ $[\nabla]$ \ of a fiber bundle over the hypersurface \ $\nabla \subset B(\C^1,3)$ \ consisting of configurations of three points \ $\alpha, \beta, \gamma \in \C^1$, \ some two of which (let us call them \ $\beta$ \ and \ $\gamma$) have equal real parts.
The fiber over such a configuration is a triangle spanned by points \ $\Phi(\alpha, \beta),$ $\Phi(\beta, \gamma)$ \ and \ $\Phi(\gamma, \alpha)$. Any such configuration defines three points in \ $\X \equiv \mbox{Sym}^2(\C^1) $, \ namely the points \ $\{\alpha, \beta\},$ $ \{\beta, \gamma\}, $ and $\{\gamma, \alpha\}$. Choosing arbitrarily some two of these three points, we obtain a point of the configuration space \ $B(\X,2)$. \ The image \ $\partial^1 ([\nabla])$ \ of the homology class of the cycle \ $[\nabla]$ \ under the homomorphism \ $\partial^1: \E^1_{1,6} \to \E^1_{0,6} \simeq \bar H_6(\Theta_0)$ \ is the homology class of the subvariety in \ $\Theta_0$ \ swept out by intervals over all points of \ $B(\X,2)$ \ obtained in this way from the points of the cycle \ $\nabla$. So, by relation \ (\ref{thomfour}) \ it remains to calculate the homology class of the set of all such points in the group \ $\bar H_5(B(\X,2), \pm \Z) \sim \Z_2$. \ The coefficient map \ $\bar H_5(B(\X,2), \pm \Z) \to \bar H_5(B(\X,2), \Z_2)$ \ is an isomorphism, therefore this class is characterized by the intersection index (mod 2) of this cycle (considered as a non-oriented one) with any compact 3-cycle generating the group \ $H_3(B(\X,2), \Z_2) \simeq \Z_2$. \ For such a 3-cycle we can take the union of pairs of opposite points of the unit sphere in \ $\X \simeq \C^2$. So, we are looking for triples of distinct points \ $\alpha, \beta, \gamma \in \C^1$ \ such that the real parts of \ $\beta$ \ and \ $\gamma$ \ are the same, and some two of three sets of coefficients $(-\alpha-\beta, \alpha \beta)$, $(-\beta-\gamma, \beta \gamma)$, $(-\gamma-\alpha, \alpha \gamma)$ \ of quadratic polynomials with roots \ $\{\alpha, \beta\}$, $\{\beta, \gamma\}$ \ and \ $\{\alpha, \gamma\}$ \ are opposite points of the unit sphere in \ $\C^2$. \ It is easy to calculate that this happens if and only if the set \ $\{\alpha, \beta, \gamma\}$ \ is equal to the set \ $\{0, i, -i\}$, \ with an arbitrary labelling of the points of the latter set by \ $\alpha, $ $\beta $ \ and \ $\gamma$. \ So, each of the three edges of the triangle over the point \ $\{0, i, -i\} \in B(\C^1,3)$ \ gives us an intersection point of our cycles. Further, \ $\{0, i, -i \}$ \ is a triple self-intersection point of the subvariety \ $\nabla \subset B(\C^1,3)$, \ therefore we have \ $3 \times 3$ \ intersections contributing to the desired intersection index, which is therefore not equal to 0. \hfill $\Box$

\section{Third term}
\label{sec3t}

Recall that the third term of the main filtration is the space of an affine bundle with base \ $\Lambda_3 \setminus \Lambda_2$ \ and fibers equal to \ $\C^{k d-3k}$. 

\begin{theorem}
\label{thmf3}
Let $\E^r_{j, \mathfrak q}$ be the spectral sequence calculating the group \ $\bar H_*(\Lambda_3 \setminus \Lambda_2)$ \
 and generated by inverse filtration \ $\{\Theta_0 \subset \Theta_1 \subset \Theta_2\}$ \ of \ $\Lambda_3 \setminus \Lambda_2$. \ Then its page
\ $\E^1$ \ is as shown in Table \ref{tab2} $($right$)$. \ Its homomorphism \ $\partial^1: \E^1_{1,11} \to \E^1_{0,11}$ \ is an isomorphism, and homomorphism \ $\partial^1: \E^1_{8,2} \to \E^1_{8,1}$ \ sends a free generator of \ $\E^1_{8,2}$ \ to an element of infinite order.
\end{theorem}

The proof of this theorem takes the rest of \S \ref{sec3t}.

\subsection{Term \ $\Theta_0$ \ of inverse filtration of \ $\Lambda_3 \setminus \Lambda_2$}
\label{t30}
The only symbol \ $A$ \ of complexity \ 3 \ and defect \ 0 \ is \ $(2,2,2)$. \ The corresponding space \ $\tilde \B(A)$ \ is the subset \
 $\tilde B(\X,3)$ \ of the configuration space \ $B(\X,3)$ \ consisting of all {\em independent} triples of elementary conditions, i.e. triples defining subspaces of complex codimension \ $3k$ \ in \ $\K$. \ Namely, it are all points of \ $B(\X,3)$ \ except for the triples of the form \ $\{\{\alpha, \beta\}, \{\beta, \gamma\}, \{\gamma, \alpha\}\}$, \ where \ $\alpha, $ $\beta$ \ and \ $\gamma$ \ are arbitrary three distinct points of \ $\C^1$. 

Term \ $\Theta_0$ \ of the inverse filtration of \ $\Lambda_3 \setminus \Lambda_2$ \ is the union of all open triangles in \ $\C^W$ \ spanned by points \ $\Phi(\chi_1), $ \ $\Phi(\chi_2)$ \ and \ $\Phi(\chi_3)$ \ for arbitrary configurations \ $\{\chi_1, \chi_2, \chi_3\} \in \tilde B(\X,3)$. \ In particular, this term is a fiber bundle with the base \ $\tilde B(\X, 3)$ \ and open triangles as fibers. The Thom isomorphism of this fiber bundle gives us the equality 
\begin{equation}
\bar H_N (\Theta_0) \simeq \bar H_{N-2}(\tilde B(\X,3), \pm \Z)
\end{equation}
 for any \ $N$. 

\begin{lemma}
\label{lem11}
The group \ $\bar H_i(B(\X,3), \pm \Z) $ \ is isomorphic to \ $\Z_2$ \ for \ $i$ \ equal to 11 and 9, to \ $\Z_3$ \ for \ $i$ \ equal to \ 10 and 6, and is trivial for all other \ $i$.
\end{lemma}

\noindent
{\it Proof.} Let us calculate groups \ $H^{12-i}(B(\R^4,3), \pm \Z)$ \ that are Poincar\'e dual to these ones. By the well-known facts about the cohomology groups (with constant coefficients) of groups \ $S(3)$ \ and \ $\Z_3$, \ and the exact sequence of the two-fold covering \ $K(\Z_3,1) \to K(S(3),1)$, \ the group \ $H^i(S(3), \pm \Z) \equiv H^i(B(\R^\infty ,3), \pm \Z)$ \ is isomorphic to \ $\Z_2$ \ for all odd \ $i$, \ to \ $\Z_3$ \ for \ $i \equiv 2 (\mbox{mod } 4),$ \ and is trivial for remaining values of \ $i$. By Theorem 2 of \cite{book}, \S I.4, the map
\begin{equation}
H^m(B(\R^\infty,3), \pm \Z) \to H^m(B(\R^4,3), \pm \Z)
\label{inclu}
\end{equation}
induced by the inclusion \ $\R^4 \hookrightarrow \R^\infty$ \ is epimorphic for all \ $i$. The canonical decomposition of the space \ $B(\R^4,3)$ \ into open cells used in this theorem has only cells of codimension up to \ $6$, \ therefore only the elements of groups \ $H^m(S(3), \pm \Z)$ \ with \ $m \leq 6$ \ can contribute to \ $H^m(B(\R^4,3), \pm \Z)$. \ All {\em stable cells} of \ $B(\R^\infty,3)$ \ of codimension up to 3 appear in \ $B(\R^4,3)$, \ therefore for \ $m \leq 3$ \ the map (\ref{inclu}) is an isomorphism. 
The only cell of codimension 6 (consisting of 3-configurations, all three points of which have equal orthogonal projections to a fixed hyperplane) appears in the boundary of each of two cells of codimension 5 with coefficient \ $\pm 3$, \ therefore the group \ $\Z_3$ \ in dimension \ $m=6$ \ survives the map (\ref{inclu}). The group \ $\Z_2$ \ in dimension 5 does not, as follows easily from the consideration of cells of codimension 4, 5 and 6. Thus, the group \ $H^m(B(\R^4,3), \pm \Z)$ \ is isomorphic to \ $\Z_2$ \ for \ $m$ \ equal to 1 and 3, to \ $\Z_3$ \ for \ $j$ \ equal to \ 2 \ and \ 6, and is trivial for all other \ $j$. \ By Poincar\'e duality, this proves our lemma. \hfill $\Box$ \medskip

The difference \ $B(\X, 3) \setminus \tilde B(\X, 3)$ \ is obviously homeomorphic to the space \ $B(\C^1, 3)$, \ and the restriction of the local system \ $\pm \Z$ \ (defined on the whole \ $B(\X, 3)$) \ to this difference is isomorphic to the local system \ $\pm \Z$ \ defined in the terms of \ $B(\C^1, 3)$. 
Therefore, the statement of Theorem \ref{thmf3} concerning the column \ $j = 0$ \ of Table \ref{tab2} follows immediately from Lemmas \ref{lem11} and \ref{lem6} and the exact sequence of Borel--Moore homology groups with coefficients in \ $\pm \Z$ \ for the pair \ $(B(\X,3), (B(\X,3) \setminus \tilde B(\X,3)))$.

\subsection{Term \ $\Theta_1 \setminus \Theta_0$ \ of inverse filtration} 
\label{term31}
This term is covered by one \ $\bar A$-block where \ $A=(3,2)$. \ It is the space of a fiber bundle, whose base \ $\tilde \B((3,2))$ \ is the space \ $\B(3,2) \equiv B(\C^1,3) \times \X$ \ from which something is removed, namely the set of pairs 
\begin{equation}
\label{32}
\left(\{\alpha, \beta, \gamma\} \in B(\C^1,3); \{\delta, \varepsilon\} \in \X \right)
\end{equation}
such that the system of conditions
\begin{equation}
\label{33}
f(\alpha)=f(\beta)=f(\gamma); \ f(\delta)=f(\varepsilon)
\end{equation}
on the map \ $f \in \K$ \ defines a subspace of codimension \ $< 3k$ \ in \ $\K$. \ It is easy to see that the last set consists of pairs (\ref{32}) 
not satisfying condition 1) of Lemma \ref{lem100}, i.e. such that both points \ $\delta \neq \varepsilon$ \ belong to the set \ $\{\alpha, \beta, \gamma\}$. \ Denote this set by \ $\triangleq $. 

For any point \ $ \Gamma \in (B(\C^1,3) \times \X) \setminus \triangleq $ \ of the form (\ref{32}), the corresponding subspace \ $L(\Gamma) \subset \K$ \ defined by conditions (\ref{33}) has codimension \ $3k$, \ and the simplex \ $\Delta(\Gamma) \subset \C^W$ \ is spanned by four points \ $\Phi(\{\alpha, \beta\}), $ $\Phi(\{\alpha, \gamma\}),$ $ \Phi(\{\beta, \gamma\})$, \ and \ $\Phi(\{\delta, \varepsilon\})$. \ The entire boundary of this simplex lies in the union of term \ $\Lambda_2$ \ of main filtration and term \ $\Theta_0$ \ of the inverse filtration of \ $\Lambda_3 \setminus \Lambda_2$. 
Thus, the space \ $\Theta_1 \setminus \Theta_0$ \ of the latter filtration is swept out by the interiors of these 3-simplices. By the Thom isomorphism, we have
\begin{equation}
\label{thom3}
\bar H_N(\Theta_1 \setminus \Theta_0) \simeq \bar H_{N-3}((B(\C^1,3) \times \X) \setminus \triangleq \ , \pm \Z),
\end{equation}
 where \ $\pm \Z$ \ is the local system on the product \ $B(\C^1,3) \times \X$ \
lifted from the eponymous local system on its first factor \ $B(\C^1,3)$. \ Let us calculate right-hand groups in (\ref{thom3}).
\medskip

By Lemma \ref{lem6} and K\"unneth formula, the group $$\bar H_i (B(\C^1,3) \times \X, \pm \Z)$$ is equal to \ $\Z_2$ \ for \ $i= 9$, \ to \ $\Z_3$ \ for \ $i= 8,$ \ and is trivial for all other \ $i$.

\begin{lemma}
\label{lem4}
1$)$ The group \ $\bar H_i(\ \triangleq \ , \pm \Z)$ \ is isomorphic to \ $\Z \oplus \Z_2$ \ for \ $i=5$, \ to \ $\Z$ \ for \ $i=4$, \ and is trivial for all other \ $i$. 

2$)$ A free generator of the group \ $\bar H_5(\ \triangleq \ , \pm \Z)$ \
 can be realized by the fundamental class of the 5-dimensional submanifold in \ $\triangleq$ \ consisting of 
pairs \ $($\ref{32}$)$ \ such that \ $\delta \neq \varepsilon$, \ $\{\delta, \varepsilon\}$ \ is a subset of \ $\{\alpha, \beta, \gamma\},$ \ and the point \ $\{\alpha, \beta, \gamma\} \setminus \{\delta, \varepsilon\}$ \ lies in the interval \ $(\delta, \varepsilon) \subset \C^1$. 
 
3$)$ For a generator of the group \ $\bar H_4(\ \triangleq \ , \pm \Z)$ \ we can take the subvariety of the previous cycle, consisting of such configurations where interval \ $(\delta, \varepsilon)$ \ is vertical. 

4$)$ The group \ $H_2(\ \triangleq \ , \pm \Z)$ \ is generated by the fundamental cycle of the submanifold in \ $\triangleq$ \ consisting of points \ $($\ref{32}$)$ \ such that \ $\alpha = \delta=0,$ \ and the points \ $ \beta = \varepsilon$ \ and \ $\gamma$ \ belong to the circles of radii \ $1$ and \ $1/2$ respectively centered at the origin in $\C^1$.

5$)$ The group \ $H_1(\ \triangleq \ , \pm \Z)$ \ is generated by the 1-cycle consisting of points \ $($\ref{32}$)$ \ such that \ $\alpha = \delta=0, $ $\beta = \varepsilon =1,$ \ and the point \ $\gamma$ \ belongs to the circle of radius \ $1/2$ \ with center at the origin in \ $\C^1$.
\end{lemma}

\noindent
{\it Proof.}
 \ $\triangleq$ \ is the space of a fiber bundle with base \ $B(\C^1,2) \sim \C^1 \times (0,\infty) \times \RP^1$, \ its fiber over any two-configuration \ $\{\delta, \varepsilon\}$ \ is equal to the space \ $\C^1 \setminus \{\delta, \varepsilon\}$ \ of choices of the point \ $\{\alpha, \beta, \gamma\} \setminus \{\delta, \varepsilon\}$. \ The rest of the calculation of its homology groups is elementary. It is also easy to check that four submanifolds indicated in statements 2)---5) as generators of corresponding homology groups are \ $\pm \Z$-orientable, and the intersection indices of these manifolds of complementary dimensions in \ $\triangleq$ \ are equal to \ $\pm 1$. \hfill $\Box$ \medskip

The structure of groups \ $\bar H_i((B(\C^1,3) \times \X)\setminus \triangleq \ , \pm \Z)$ \ now follows from the exact sequence of the pair \ $((B(\C^1,3) \times \X), \triangleq \ )$. \ Namely, these groups are equal to \ $\Z_2$ \ for \ $i= 9$, \ to \ $\Z_3$ \ for \ $i=8$, \ to \ $\Z_2 \oplus \Z$ \ for \ $i= 6$, \ and to \ $\Z$ \ for \ $i= 5$. \ Assertion of Theorem \ref{thmf3} about the column \ $ j= 1$ \ follows immediately from this calculation and identity (\ref{thom3}). \medskip

The statement of Theorem \ref{thmf3} about the differential \ $\partial^1: \E^1_{1,11} \to \E^1_{0,11}$ \ can be proved in the same way as the statement of Theorem \ref{thmf2} about the map \ $\partial^1: \E^1_{1,6} \to \E^1_{0,6}$.

\subsection{Term \ $\Theta_2 \setminus \Theta_1$ \ of inverse filtration}

The only symbol \ $A$ \ of complexity \ 3 \ and defect \ 2 \ is \ $(4)$, \ so the part of \ $\Lambda_3 \setminus \Lambda_2$ \ not covered by the above described subset \ $\Theta_1$ \ is covered by the \ $(4)$-block. This block is the space of a fiber bundle over the space \ $B(\C^1,4)$ \ of quadruples of distinct points in \ $\C^1.$ \ Its fiber over such a point \ $\Gamma = \{\alpha, \beta, \gamma, \delta\} \subset \C^1$ \ is the 5-dimensional simplex \ $\Delta(\Gamma)$, \ six vertices \ $\Phi(\{\alpha, \beta\}),$ \ $ \Phi(\{\alpha, \gamma\})$, \ etc. of which correspond to two-element subsets of the set \ $\{\alpha, \beta, \gamma, \delta\}$. \
According to Lemma \ref{lem18}, the part of such a simplex \ $\Delta(\Gamma)$ \ that lies in \ $\Theta_2 \setminus \Theta_1$ \ is the union of the interior points of its faces corresponding to two-connected simple graphs on four vertices \ $\alpha, \beta, \gamma$ \ and \ $\delta$. Namely, these are all graphs with six or five edges, and all three graphs with four edges that define a 4-cycle. Thus, term \ $\Theta_2 \setminus \Theta_1$ \ is the space of a fiber bundle over \ $B(\C^1,4)$, \ whose fiber over the configuration \ $\{\alpha, \beta, \gamma, \delta\} \in B(\C^1,4)$ \ is the union of interior points of all faces of the corresponding 5-simplex, which correspond to two-connected graphs on four vertices. 
 An elementary calculation (see \cite{twocon}) gives us the following specialization of Theorem \ref{2c}.

\begin{lemma}
\label{2congrlem}
The Borel--Moore homology group of such a fiber is non-trivial only in dimension \ 4 \ and is isomorphic to \ $\Z^2:$ \ it is generated by the classes of arbitrary two of the three chains 
\unitlength 0.6mm
\begin{equation}
\label{2conhom}
\mbox{\begin{picture}(30,13)
\put(0,0){$\beta$}
\put(0,10){$\alpha$}
\put(16,0){$\gamma$}
\put(16,10){$\delta$}
\put(5,0){\circle*{1.5}}
\put(5,10){\circle*{1.5}}
\put(15,0){\circle*{1.5}}
\put(15,10){\circle*{1.5}}
\put(5,0){\line(1,0){10}}
\put(5,0){\line(1,1){10}}
\put(5,10){\line(1,0){10}}
\put(15,0){\line(0,1){10}}
\put(15,0){\line(-1,1){10}}
\put(22,5){$-$}
\end{picture}} 
\mbox{\begin{picture}(20,13)
\put(0,0){$\beta$}
\put(0,10){$\alpha$}
\put(16,0){$\gamma$}
\put(16,10){$\delta$}
\put(5,0){\circle*{1.5}}
\put(5,10){\circle*{1.5}}
\put(15,0){\circle*{1.5}}
\put(15,10){\circle*{1.5}}
\put(5,0){\line(1,0){10}}
\put(5,0){\line(1,1){10}}
\put(5,0){\line(0,1){10}}
\put(5,10){\line(1,0){10}}
\put(15,0){\line(-1,1){10}}
\end{picture}} \ , \qquad
\mbox{\begin{picture}(30,13)
\put(0,0){$\beta$}
\put(0,10){$\alpha$}
\put(16,0){$\gamma$}
\put(16,10){$\delta$}
\put(5,0){\circle*{1.5}}
\put(5,10){\circle*{1.5}}
\put(15,0){\circle*{1.5}}
\put(15,10){\circle*{1.5}}
\put(5,0){\line(1,0){10}}
\put(5,0){\line(0,1){10}}
\put(5,10){\line(1,0){10}}
\put(15,0){\line(0,1){10}}
\put(15,0){\line(-1,1){10}}
\put(22,5){$-$}
\end{picture}} 
\mbox{\begin{picture}(20,13)
\put(0,0){$\beta$}
\put(0,10){$\alpha$}
\put(16,0){$\gamma$}
\put(16,10){$\delta$}
\put(5,0){\circle*{1.5}}
\put(5,10){\circle*{1.5}}
\put(15,0){\circle*{1.5}}
\put(15,10){\circle*{1.5}}
\put(5,0){\line(1,0){10}}
\put(5,0){\line(1,1){10}}
\put(5,0){\line(0,1){10}}
\put(5,10){\line(1,0){10}}
\put(15,0){\line(0,1){10}}
\end{picture}} \ , \qquad
 \mbox{\begin{picture}(30,13)
\put(0,0){$\beta$}
\put(0,10){$\alpha$}
\put(16,0){$\gamma$}
\put(16,10){$\delta$}
\put(5,0){\circle*{1.5}}
\put(5,10){\circle*{1.5}}
\put(15,0){\circle*{1.5}}
\put(15,10){\circle*{1.5}}
\put(5,0){\line(1,1){10}}
\put(5,0){\line(1,0){10}}
\put(5,0){\line(0,1){10}}
\put(15,0){\line(0,1){10}}
\put(15,0){\line(-1,1){10}}
\put(22,5){$-$}
\end{picture}}
\mbox{\begin{picture}(20,13)
\put(0,0){$\beta$}
\put(0,10){$\alpha$}
\put(16,0){$\gamma$}
\put(16,10){$\delta$}
\put(5,0){\circle*{1.5}}
\put(5,10){\circle*{1.5}}
\put(15,0){\circle*{1.5}}
\put(15,10){\circle*{1.5}}
\put(5,0){\line(1,1){10}}
\put(5,0){\line(0,1){10}}
\put(5,10){\line(1,0){10}}
\put(15,0){\line(0,1){10}}
\put(15,0){\line(-1,1){10}}
\end{picture}}
\end{equation}
$($where the orientations of the faces represented by these graphs are determined by the order \ $\alpha < \beta < \gamma < \delta)$. \ The sum of all these three chains is equal to the boundary of entire simplex $($depicted by the complete graph$)$. \hfill $\Box$
\label{lem22}
\end{lemma}

\begin{remark} \rm
\label{ir}
There is an important one-to-one correspondence between these three basic chains and matchings in the set \ $\{\alpha, \beta, \gamma, \delta \}$: \ say, the chain \ (\ref{2conhom}) \ containing graphs with missing edge \ $(\alpha, \gamma)$ \ or \ $(\beta, \delta)$ \ corresponds to the matching \ $(\alpha, \gamma)(\beta, \delta)$.
\end{remark}

Let \ $A_2$ \ be the quotient lattice of \ $\Z^3$ \ along the diagonal sublattice consisting of the points \ $(t, t, t)$. 
Denote by \ $\{A_2\}$ \ the representation of the braid group \ $Br_4 \equiv \pi_1(B(\C^1,4))$ \ in \ $A_2$ \ defined as the composition
\begin{equation}
\label{compos}
Br_4 \to S(4) \to S(3) \to \mbox{Aut}(A_2) ,
\end{equation}
where the first homomorphism takes any braid to the corresponding permutation of points of a distinguished \ 4-configuration, the second one takes any rotation of a tetrahedron to the corresponding permutation of pairs of its opposite edges, and the last one is the standard action of \ $S(3)$ \ in $A_2$ \ by permutations of coordinates in \ $\Z^3$. 

\begin{lemma}
\label{lem13}  
$\bar H_N(\Theta_2 \setminus \Theta_1) \simeq \bar H_{N-4}(B(\C^1,4), \{A_2\}) $ for any natural \ $N$.
\end{lemma}

\noindent
{\it Proof.} It is easy to check that any permutation of points \ $\alpha, \beta, \gamma, \delta$ \ acts on the set of cycles (\ref{2conhom}) in exactly the same way as
the automorphism of \ $A_2$ \ obtained from this permutation by the composition of two last arrows in (\ref{compos}) acts on the images of the basic vectors of \ $\Z^3$ \ under the factorization by the diagonal. Namely, any transposition of some two of these four points preserves the cycle containing the graph with missing edge connecting these points, and permutes two other generators. \hfill $\Box$ \medskip

Let also \ $\widehat{\{A_2\}}$ \ be a representation of \ $Br_4$ \ in \ $\mbox{Aut}(\Z^3)$ \ defined by a composition analogous to (\ref{compos}), in which the last homomorphism is the action by permutations on the whole \ $\Z^3$ \ and not on its quotient lattice \ $A_2$.

\begin{theorem}
\label{homcal}
All non-trivial groups \ $\bar H_i(B(\C^1,4), \Z)$, \ $\bar H_i(B(\C^1,4), \widehat{\{A_2\}})$ \ and \ $\bar H_{i}(B(\C^1,4), \{A_2\})$ \ are as shown in Table \ref{tab5}.
\end{theorem}

The column \ $j = 2$ \ of Table \ref{tab2} (right) follows immediately from Lemma \ref{lem13} and the last row of Table \ref{tab5}.

\begin{table}
\caption{}
\label{tab5}
\begin{tabular}{|c|c|c|c|c|}
\hline
$i$ & 8 & 7 &6 & 5 \\
\hline
$\bar H_{i}(B(\C^1,4), \Z)$ \ \phantom{\LARGE $A$} \ & $\Z$ & $\Z$ & 0 & $\Z_2$ \\
\hline
 $\bar H_i(B(\C^1,4), \widehat{\{A_2\}})$ \phantom{\LARGE $I$} & $\Z$ & $\Z^2$ & $\Z \oplus \Z_2$ & $(\Z_2)^2$ \\
\hline
$\bar H_{i}(B(\C^1,4), \{A_2\})$ \phantom{\LARGE $I$}& 0 & $\Z$ & $\Z \oplus \Z_2$ & $\Z_2$ \\
\hline
\end{tabular}

\end{table}

\medskip

\noindent
{\it Proof of Theorem \ref{homcal}}. 
The second row of Table \ref{tab5} is proved in \cite{AA}.

Consider the decomposition of the space \ $B(\C^1,4)$ \ into open cells used in \cite{fuks}. Namely, we denote by \
 \unitlength 0.5mm
\begin{picture}(21,10)
\put(10.5,3){\oval(21,10)}
\put(3,3){\circle*{2}} 
\put(8,3){\circle*{2}} 
\put(13,3){\circle*{2}} 
\put(18,3){\circle*{2}} 
\end{picture}
\ the cell consisting of configurations, all four points of which have different real parts, and by \
\begin{picture}(16,10)
\put(8,3){\oval(16,10)}
\put(3,1){\circle*{2}} 
\put(3,5){\circle*{2}} 
\put(8,3){\circle*{2}} 
\put(13,3){\circle*{2}} 
\end{picture} 
\ the 7-dimensional cell (called \ $e(2,1,1)$ \ in \cite{fuks}) consisting of configurations in which only the two leftmost points have equal real parts. Analogously, denote by \
\begin{picture}(16,10)
\put(8,3){\oval(16,10)}
\put(3,3){\circle*{2}} 
\put(8,1){\circle*{2}} 
\put(8,5){\circle*{2}} 
\put(13,3){\circle*{2}} 
\end{picture}, 
\begin{picture}(16,10)
\put(8,3){\oval(16,10)}
\put(3,3){\circle*{2}} 
\put(8,3){\circle*{2}} 
\put(13,1){\circle*{2}} 
\put(13,5){\circle*{2}} 
\end{picture}, 
\begin{picture}(12,10)
\put(6,3){\oval(12,10)}
\put(4,0){\circle*{2}} 
\put(4,3){\circle*{2}} 
\put(4,6){\circle*{2}} 
\put(9,3){\circle*{2}} 
\end{picture}, 
\begin{picture}(11,10)
\put(5.5,3){\oval(11,10)}
\put(3,1){\circle*{2}} 
\put(3,4){\circle*{2}} 
\put(8,2){\circle*{2}} 
\put(8,5){\circle*{2}} 
\end{picture}, 
\begin{picture}(12,10)
\put(6,3){\oval(12,10)}
\put(3,3){\circle*{2}} 
\put(8,0){\circle*{2}} 
\put(8,3){\circle*{2}} 
\put(8,6){\circle*{2}} 
\end{picture}, \ and \
\begin{picture}(8,10)
\put(4,4){\oval(8,13)}
\put(4,-1){\circle*{2}} 
\put(4,2){\circle*{2}} 
\put(4,5){\circle*{2}} 
\put(4,8){\circle*{2}} 
\end{picture}
\ the cells called in \cite{fuks} respectively \
 $e(1,2,1), $ $e(1,1,2)$, $e(3,1), $ $e(2,2)$, $e(1,3)$, \ and \ $e(4)$. \ In particular, \ \begin{picture}(8,10)
\put(4,4){\oval(8,13)}
\put(4,-1){\circle*{2}} 
\put(4,2){\circle*{2}} 
\put(4,5){\circle*{2}} 
\put(4,8){\circle*{2}} 
\end{picture} \ is the only 5-dimensional cell of this decomposition: it consists of configurations, all whose points have equal real parts. 

$\bar H_*(B(\C^1,4), \widehat{\{A_2\}})$ \ is the homology group of a 3-fold covering over \ $B(\C^1,4)$. The decomposition of this covering space into cells is lifted from the previous one: the preimage of each of our cells \
\unitlength 0.5mm
\begin{picture}(21,10)
\put(10.5,3){\oval(21,10)}
\put(3,3){\circle*{2}} 
\put(8,3){\circle*{2}} 
\put(13,3){\circle*{2}} 
\put(18,3){\circle*{2}} 
\end{picture}, \
\begin{picture}(16,10)
\put(8,3){\oval(16,10)}
\put(3,1){\circle*{2}} 
\put(3,5){\circle*{2}} 
\put(8,3){\circle*{2}} 
\put(13,3){\circle*{2}} 
\end{picture}, \
 \begin{picture}(16,10)
\put(8,3){\oval(16,10)}
\put(3,3){\circle*{2}} 
\put(8,1){\circle*{2}} 
\put(8,5){\circle*{2}} 
\put(13,3){\circle*{2}} 
\end{picture} , \dots , \ \begin{picture}(8,10)
\put(4,4){\oval(8,13)}
\put(4,-1){\circle*{2}} 
\put(4,2){\circle*{2}} 
\put(4,5){\circle*{2}} 
\put(4,8){\circle*{2}} 
\end{picture} \ 
consists of three cells, the notation of which is obtained from that of the original cell by adding a subscript equal to 2, 3 or 4. Namely, these three cells are in one-to-one correspondence with matchings of points of an arbitrary $4$-configuration \ $\{\alpha, \beta, \gamma, \delta\}$ \ from the original cell. We order lexicographically the points of any such 4-configuration: 
 first by increase of their real parts and then by decrease of imaginary parts. This order is continuous along any of our cells. Then, given a cell augmented by a matching of four points of corresponding 4-configurations, we mark it with the subscript equal to the number of the point matched with the first one.

\begin{definition} \rm
\label{rem1}
The {\em standard orientation} of any our cells is determined by our 
 lexicographic order of points of configurations from this cell and is given by the following differential form:
$$\begin{array}{l}
d(\mbox{the smallest real value of its points}) \wedge \\
\wedge d(\mbox{the next smallest real value}) \wedge \dots 
\wedge d(\mbox{the largest real value}) \wedge \\
\wedge d(\mbox{the imaginary part of the first ordered point}) \wedge \dots \wedge \\
\wedge d(\mbox{the imaginary part of the fourth point}) ,
\end{array}
$$
 cf. \cite{vai}.
\end{definition}

It is easy to calculate that differentials in the arising complex are as follows: 

$\partial \left(\mbox{
\begin{picture}(21,10)
\put(10.5,3){\oval(21,10)}
\put(3,3){\circle*{2}} 
\put(8,3){\circle*{2}} 
 \put(13,3){\circle*{2}} 
\put(18,3){\circle*{2}} 
\end{picture}}_2 \right) = - 
\mbox{\begin{picture}(16,10)
\put(8,3){\oval(16,10)}
\put(3,3){\circle*{2}} 
\put(8,1){\circle*{2}} 
\put(8,5){\circle*{2}} 
\put(13,3){\circle*{2}} 
\end{picture}}_2+
\mbox{\begin{picture}(16,10)
\put(8,3){\oval(16,10)}
\put(3,3){\circle*{2}} 
\put(8,1){\circle*{2}} 
\put(8,5){\circle*{2}} 
\put(13,3){\circle*{2}} 
\end{picture}}_3$

$\partial \left(\mbox{
\begin{picture}(21,10)
\put(10.5,3){\oval(21,10)}
\put(3,3){\circle*{2}} 
\put(8,3){\circle*{2}} 
\put(13,3){\circle*{2}} 
\put(18,3){\circle*{2}} 
\end{picture}}_3 \right) = 
\mbox{\begin{picture}(16,10)
\put(8,3){\oval(16,10)}
\put(3,1){\circle*{2}} 
\put(3,5){\circle*{2}} 
\put(8,3){\circle*{2}} 
\put(13,3){\circle*{2}} 
\end{picture}}_3 - 
\mbox{\begin{picture}(16,10)
\put(8,3){\oval(16,10)}
\put(3,1){\circle*{2}} 
\put(3,5){\circle*{2}} 
\put(8,3){\circle*{2}} 
\put(13,3){\circle*{2}} 
\end{picture}}_4- \mbox{\begin{picture}(16,10)
\put(8,3){\oval(16,10)}
\put(3,3){\circle*{2}} 
\put(8,1){\circle*{2}} 
\put(8,5){\circle*{2}} 
\put(13,3){\circle*{2}} 
\end{picture}}_3+ \mbox{\begin{picture}(16,10)
\put(8,3){\oval(16,10)}
\put(3,3){\circle*{2}} 
\put(8,1){\circle*{2}} 
\put(8,5){\circle*{2}} 
\put(13,3){\circle*{2}} 
\end{picture}}_2 + \mbox{\begin{picture}(16,10)
\put(8,3){\oval(16,10)}
\put(3,3){\circle*{2}} 
\put(8,3){\circle*{2}} 
\put(13,1){\circle*{2}} 
\put(13,5){\circle*{2}} 
\end{picture}}_3-\mbox{\begin{picture}(16,10)
\put(8,3){\oval(16,10)}
\put(3,3){\circle*{2}} 
\put(8,3){\circle*{2}} 
\put(13,1){\circle*{2}} 
\put(13,5){\circle*{2}} 
\end{picture}}_4 $

$\partial \left(\mbox{
\begin{picture}(21,10)
\put(10.5,3){\oval(21,10)}
\put(3,3){\circle*{2}} 
\put(8,3){\circle*{2}} 
\put(13,3){\circle*{2}} 
\put(18,3){\circle*{2}} 
\end{picture}}_4 \right) = 
 \mbox{\begin{picture}(16,10)
\put(8,3){\oval(16,10)}
\put(3,1){\circle*{2}} 
\put(3,5){\circle*{2}} 
\put(8,3){\circle*{2}} 
\put(13,3){\circle*{2}} 
\end{picture}}_4 - 
\mbox{\begin{picture}(16,10)
\put(8,3){\oval(16,10)}
\put(3,1){\circle*{2}} 
\put(3,5){\circle*{2}} 
\put(8,3){\circle*{2}} 
\put(13,3){\circle*{2}} 
\end{picture}}_3 + \mbox{\begin{picture}(16,10)
\put(8,3){\oval(16,10)}
\put(3,3){\circle*{2}} 
\put(8,3){\circle*{2}} 
\put(13,1){\circle*{2}} 
\put(13,5){\circle*{2}} 
\end{picture}}_4 - \mbox{\begin{picture}(16,10)
\put(8,3){\oval(16,10)}
\put(3,3){\circle*{2}} 
\put(8,3){\circle*{2}} 
\put(13,1){\circle*{2}} 
\put(13,5){\circle*{2}} 
\end{picture}}_3$

$\partial \left( \mbox{\begin{picture}(16,10)
\put(8,3){\oval(16,10)}
\put(3,1){\circle*{2}} 
\put(3,5){\circle*{2}} 
\put(8,3){\circle*{2}} 
\put(13,3){\circle*{2}} 
\end{picture}}_2 \right) =\ 
\mbox{\begin{picture}(12,10)
\put(6,3){\oval(12,10)}
\put(4,0){\circle*{2}} 
\put(4,3){\circle*{2}} 
\put(4,6){\circle*{2}} 
\put(9,3){\circle*{2}} 
\end{picture}}_2-\ \mbox{\begin{picture}(12,10)
\put(6,3){\oval(12,10)}
\put(4,0){\circle*{2}} 
\put(4,3){\circle*{2}} 
\put(4,6){\circle*{2}} 
\put(9,3){\circle*{2}} 
\end{picture}}_3+ \mbox{\begin{picture}(12,10)
\put(6,3){\oval(12,10)}
\put(4,0){\circle*{2}} 
\put(4,3){\circle*{2}} 
\put(4,6){\circle*{2}} 
\put(9,3){\circle*{2}} 
\end{picture}}_4 $

$\partial \left( \mbox{\begin{picture}(16,10)
\put(8,3){\oval(16,10)}
\put(3,1){\circle*{2}} 
\put(3,5){\circle*{2}} 
\put(8,3){\circle*{2}} 
\put(13,3){\circle*{2}} 
\end{picture}}_3 \right) = \mbox{\begin{picture}(12,10)
\put(6,3){\oval(12,10)}
\put(4,0){\circle*{2}} 
\put(4,3){\circle*{2}} 
\put(4,6){\circle*{2}} 
\put(9,3){\circle*{2}} 
\end{picture}}_3-\ \mbox{\begin{picture}(11,10)
\put(5.5,3){\oval(12,10)}
\put(3,1){\circle*{2}} 
\put(3,4){\circle*{2}} 
\put(8,2){\circle*{2}} 
\put(8,5){\circle*{2}} 
\end{picture}}_3+\ \mbox{\begin{picture}(11,10)
\put(5.5,3){\oval(12,10)}
\put(3,1){\circle*{2}} 
\put(3,4){\circle*{2}} 
\put(8,2){\circle*{2}} 
\put(8,5){\circle*{2}} 
\end{picture}}_4$

$\partial \left( \mbox{\begin{picture}(16,10)
\put(8,3){\oval(16,10)}
\put(3,1){\circle*{2}} 
\put(3,5){\circle*{2}} 
\put(8,3){\circle*{2}} 
\put(13,3){\circle*{2}} 
\end{picture}}_4 \right) = \ \mbox{\begin{picture}(12,10)
\put(6,3){\oval(12,10)}
\put(4,0){\circle*{2}} 
\put(4,3){\circle*{2}} 
\put(4,6){\circle*{2}} 
\put(9,3){\circle*{2}} 
\end{picture}}_3-\ \mbox{\begin{picture}(11,10)
\put(5.5,3){\oval(12,10)}
\put(3,1){\circle*{2}} 
\put(3,4){\circle*{2}} 
\put(8,2){\circle*{2}} 
\put(8,5){\circle*{2}} 
\end{picture}}_4+\ \mbox{\begin{picture}(11,10)
\put(5.5,3){\oval(12,10)}
\put(3,1){\circle*{2}} 
\put(3,4){\circle*{2}} 
\put(8,2){\circle*{2}} 
\put(8,5){\circle*{2}} 
\end{picture}}_3 $

$\partial \left(\mbox{\begin{picture}(16,10)
\put(8,3){\oval(16,10)}
\put(3,3){\circle*{2}} 
\put(8,1){\circle*{2}} 
\put(8,5){\circle*{2}} 
\put(13,3){\circle*{2}} 
\end{picture}}_2 \right) = \ \mbox{\begin{picture}(12,10)
\put(6,3){\oval(12,10)}
\put(4,0){\circle*{2}} 
\put(4,3){\circle*{2}} 
\put(4,6){\circle*{2}} 
\put(9,3){\circle*{2}} 
\end{picture}}_3-\ \mbox{\begin{picture}(12,10)
\put(6,3){\oval(12,10)}
\put(3,3){\circle*{2}} 
\put(8,0){\circle*{2}} 
\put(8,3){\circle*{2}} 
\put(8,6){\circle*{2}} 
\end{picture}}_3$

$\partial \left( \mbox{\begin{picture}(16,10)
\put(8,3){\oval(16,10)}
\put(3,3){\circle*{2}} 
\put(8,1){\circle*{2}} 
\put(8,5){\circle*{2}} 
\put(13,3){\circle*{2}} 
\end{picture}}_3 \right) = \ \mbox{\begin{picture}(12,10)
\put(6,3){\oval(12,10)}
\put(4,0){\circle*{2}} 
\put(4,3){\circle*{2}} 
\put(4,6){\circle*{2}} 
\put(9,3){\circle*{2}} 
\end{picture}}_3-\ \mbox{\begin{picture}(12,10)
\put(6,3){\oval(12,10)}
\put(3,3){\circle*{2}} 
\put(8,0){\circle*{2}} 
\put(8,3){\circle*{2}} 
\put(8,6){\circle*{2}} 
\end{picture}}_3$

$\partial \left( \mbox{\begin{picture}(16,10)
\put(8,3){\oval(16,10)}
\put(3,3){\circle*{2}} 
\put(8,1){\circle*{2}} 
\put(8,5){\circle*{2}} 
\put(13,3){\circle*{2}} 
\end{picture}}_4 \right) = \ \mbox{\begin{picture}(12,10)
\put(6,3){\oval(12,10)}
\put(4,0){\circle*{2}} 
\put(4,3){\circle*{2}} 
\put(4,6){\circle*{2}} 
\put(9,3){\circle*{2}} 
\end{picture}}_4-\ \mbox{\begin{picture}(12,10)
\put(6,3){\oval(12,10)}
\put(4,0){\circle*{2}} 
\put(4,3){\circle*{2}} 
\put(4,6){\circle*{2}} 
\put(9,3){\circle*{2}} 
\end{picture}}_3+\ \mbox{\begin{picture}(12,10)
\put(6,3){\oval(12,10)}
\put(4,0){\circle*{2}} 
\put(4,3){\circle*{2}} 
\put(4,6){\circle*{2}} 
\put(9,3){\circle*{2}} 
\end{picture}}_2-\ \mbox{\begin{picture}(12,10)
\put(6,3){\oval(12,10)}
\put(3,3){\circle*{2}} 
\put(8,0){\circle*{2}} 
\put(8,3){\circle*{2}} 
\put(8,6){\circle*{2}} 
\end{picture}}_4+\ \mbox{\begin{picture}(12,10)
\put(6,3){\oval(12,10)}
\put(3,3){\circle*{2}} 
\put(8,0){\circle*{2}} 
\put(8,3){\circle*{2}} 
\put(8,6){\circle*{2}} 
\end{picture}}_3-\ \mbox{\begin{picture}(12,10)
\put(6,3){\oval(12,10)}
\put(3,3){\circle*{2}} 
\put(8,0){\circle*{2}} 
\put(8,3){\circle*{2}} 
\put(8,6){\circle*{2}} 
\end{picture}}_2$

$\partial \left(\mbox{\begin{picture}(16,10)
\put(8,3){\oval(16,10)}
\put(3,3){\circle*{2}} 
\put(8,3){\circle*{2}} 
\put(13,1){\circle*{2}} 
\put(13,5){\circle*{2}} 
\end{picture}}_2 \right) = -\ \mbox{\begin{picture}(12,10)
\put(6,3){\oval(12,10)}
\put(3,3){\circle*{2}} 
\put(8,0){\circle*{2}} 
\put(8,3){\circle*{2}} 
\put(8,6){\circle*{2}} 
\end{picture}}_2+\ \mbox{\begin{picture}(12,10)
\put(6,3){\oval(12,10)}
\put(3,3){\circle*{2}} 
\put(8,0){\circle*{2}} 
\put(8,3){\circle*{2}} 
\put(8,6){\circle*{2}} 
\end{picture}}_3-\ \mbox{\begin{picture}(12,10)
\put(6,3){\oval(12,10)}
\put(3,3){\circle*{2}} 
\put(8,0){\circle*{2}} 
\put(8,3){\circle*{2}} 
\put(8,6){\circle*{2}} 
\end{picture}}_4 $

$\partial \left(\mbox{\begin{picture}(16,10)
\put(8,3){\oval(16,10)}
\put(3,3){\circle*{2}} 
\put(8,3){\circle*{2}} 
\put(13,1){\circle*{2}} 
\put(13,5){\circle*{2}} 
\end{picture}}_3 \right) = \ \mbox{\begin{picture}(11,10)
\put(5.5,3){\oval(11,10)}
\put(3,1){\circle*{2}} 
\put(3,4){\circle*{2}} 
\put(8,2){\circle*{2}} 
\put(8,5){\circle*{2}} 
\end{picture}}_3-\ \mbox{\begin{picture}(11,10)
\put(5.5,3){\oval(11,10)}
\put(3,1){\circle*{2}} 
\put(3,4){\circle*{2}} 
\put(8,2){\circle*{2}} 
\put(8,5){\circle*{2}} 
\end{picture}}_4-\ \mbox{\begin{picture}(12,10)
\put(6,3){\oval(12,10)}
\put(3,3){\circle*{2}} 
\put(8,0){\circle*{2}} 
\put(8,3){\circle*{2}} 
\put(8,6){\circle*{2}} 
\end{picture}}_3$

$\partial \left(\mbox{\begin{picture}(16,10)
\put(8,3){\oval(16,10)}
\put(3,3){\circle*{2}} 
\put(8,3){\circle*{2}} 
\put(13,1){\circle*{2}} 
\put(13,5){\circle*{2}} 
\end{picture}}_4 \right) = \ \mbox{\begin{picture}(11,10)
\put(5.5,3){\oval(11,10)}
\put(3,1){\circle*{2}} 
\put(3,4){\circle*{2}} 
\put(8,2){\circle*{2}} 
\put(8,5){\circle*{2}} 
\end{picture}}_4-\ \mbox{\begin{picture}(11,10)
\put(5.5,3){\oval(11,10)}
\put(3,1){\circle*{2}} 
\put(3,4){\circle*{2}} 
\put(8,2){\circle*{2}} 
\put(8,5){\circle*{2}} 
\end{picture}}_3-\ \mbox{\begin{picture}(12,10)
\put(6,3){\oval(12,10)}
\put(3,3){\circle*{2}} 
\put(8,0){\circle*{2}} 
\put(8,3){\circle*{2}} 
\put(8,6){\circle*{2}} 
\end{picture}}_3 $

$\partial \left( \mbox{\begin{picture}(12,10)
\put(6,3){\oval(12,10)}
\put(4,0){\circle*{2}} 
\put(4,3){\circle*{2}} 
\put(4,6){\circle*{2}} 
\put(9,3){\circle*{2}} 
\end{picture}}_2 \right) = \mbox{\begin{picture}(8,10)
\put(4,4){\oval(8,13)}
\put(4,-1){\circle*{2}} 
\put(4,2){\circle*{2}} 
\put(4,5){\circle*{2}} 
\put(4,8){\circle*{2}} 
\end{picture}}_3 - \mbox{\begin{picture}(8,10)
\put(4,4){\oval(8,13)}
\put(4,-1){\circle*{2}} 
\put(4,2){\circle*{2}} 
\put(4,5){\circle*{2}} 
\put(4,8){\circle*{2}} 
\end{picture}}_4$

$\partial \left( \mbox{\begin{picture}(12,10)
\put(6,3){\oval(12,10)}
\put(4,0){\circle*{2}} 
\put(4,3){\circle*{2}} 
\put(4,6){\circle*{2}} 
\put(9,3){\circle*{2}} 
\end{picture}}_3\right) = 0 $

$\partial \left( \mbox{\begin{picture}(12,10)
\put(6,3){\oval(12,10)}
\put(4,0){\circle*{2}} 
\put(4,3){\circle*{2}} 
\put(4,6){\circle*{2}} 
\put(9,3){\circle*{2}} 
\end{picture}}_4\right) = \mbox{\begin{picture}(8,10)
\put(4,4){\oval(8,13)}
\put(4,-1){\circle*{2}} 
\put(4,2){\circle*{2}} 
\put(4,5){\circle*{2}} 
\put(4,8){\circle*{2}} 
\end{picture}}_4 - \mbox{\begin{picture}(8,10)
\put(4,4){\oval(8,13)}
\put(4,-1){\circle*{2}} 
\put(4,2){\circle*{2}} 
\put(4,5){\circle*{2}} 
\put(4,8){\circle*{2}} 
\end{picture}}_3 $

$\partial \left( \mbox{\begin{picture}(11,10)
\put(5.5,3){\oval(11,10)}
\put(3,1){\circle*{2}} 
\put(3,4){\circle*{2}} 
\put(8,2){\circle*{2}} 
\put(8,5){\circle*{2}} 
\end{picture}}_2 \right) = 2 \left(\mbox{\begin{picture}(8,12)
\put(4,4){\oval(8,13)}
\put(4,-1){\circle*{2}} 
\put(4,2){\circle*{2}} 
\put(4,5){\circle*{2}} 
\put(4,8){\circle*{2}} 
\end{picture}}_2 - \mbox{\begin{picture}(8,12)
\put(4,4){\oval(8,13)}
\put(4,-1){\circle*{2}} 
\put(4,2){\circle*{2}} 
\put(4,5){\circle*{2}} 
\put(4,8){\circle*{2}} 
\end{picture}}_3 + \mbox{\begin{picture}(8,12)
\put(4,4){\oval(8,13)}
\put(4,-1){\circle*{2}} 
\put(4,2){\circle*{2}} 
\put(4,5){\circle*{2}} 
\put(4,8){\circle*{2}} 
\end{picture}}_4\right)$

$\partial \left( \mbox{\begin{picture}(11,10)
\put(5.5,3){\oval(11,10)}
\put(3,1){\circle*{2}} 
\put(3,4){\circle*{2}} 
\put(8,2){\circle*{2}} 
\put(8,5){\circle*{2}} 
\end{picture}}_3 \right) = 2 \ \mbox{\begin{picture}(8,12)
\put(4,4){\oval(8,13)}
\put(4,-1){\circle*{2}} 
\put(4,2){\circle*{2}} 
\put(4,5){\circle*{2}} 
\put(4,8){\circle*{2}} 
\end{picture}}_3$

$\partial \left( \mbox{\begin{picture}(11,10)
\put(5.5,3){\oval(11,10)}
\put(3,1){\circle*{2}} 
\put(3,4){\circle*{2}} 
\put(8,2){\circle*{2}} 
\put(8,5){\circle*{2}} 
\end{picture}}_4 \right) = 2 \ \mbox{\begin{picture}(8,12)
\put(4,4){\oval(8,13)}
\put(4,-1){\circle*{2}} 
\put(4,2){\circle*{2}} 
\put(4,5){\circle*{2}} 
\put(4,8){\circle*{2}} 
\end{picture}}_3$

$\partial \left( \mbox{\begin{picture}(12,10)
\put(6,3){\oval(12,10)}
\put(3,3){\circle*{2}} 
\put(8,0){\circle*{2}} 
\put(8,3){\circle*{2}} 
\put(8,6){\circle*{2}} 
\end{picture}}_2 \right) = \mbox{\begin{picture}(8,10)
\put(4,4){\oval(8,13)}
\put(4,-1){\circle*{2}} 
\put(4,2){\circle*{2}} 
\put(4,5){\circle*{2}} 
\put(4,8){\circle*{2}} 
\end{picture}}_3 - \mbox{\begin{picture}(8,10)
\put(4,4){\oval(8,13)}
\put(4,-1){\circle*{2}} 
\put(4,2){\circle*{2}} 
\put(4,5){\circle*{2}} 
\put(4,8){\circle*{2}} 
\end{picture}}_4$

$\partial \left( \mbox{\begin{picture}(12,10)
\put(6,3){\oval(12,10)}
\put(3,3){\circle*{2}} 
\put(8,0){\circle*{2}} 
\put(8,3){\circle*{2}} 
\put(8,6){\circle*{2}} 
\end{picture}}_3 \right) = 0$

$\partial \left( \mbox{\begin{picture}(12,10)
\put(6,3){\oval(12,10)}
\put(3,3){\circle*{2}} 
\put(8,0){\circle*{2}} 
\put(8,3){\circle*{2}} 
\put(8,6){\circle*{2}} 
\end{picture}}_4 \right) = \mbox{\begin{picture}(8,10)
\put(4,4){\oval(8,13)}
\put(4,-1){\circle*{2}} 
\put(4,2){\circle*{2}} 
\put(4,5){\circle*{2}} 
\put(4,8){\circle*{2}} 
\end{picture}}_4 - \mbox{\begin{picture}(8,10)
\put(4,4){\oval(8,13)}
\put(4,-1){\circle*{2}} 
\put(4,2){\circle*{2}} 
\put(4,5){\circle*{2}} 
\put(4,8){\circle*{2}} 
\end{picture}}_3$
\medskip

The statement of Theorem \ref{homcal} concerning groups \ $\bar H_i(B(\C^1,4), \widehat{\{A_2\}})$ \ follows immediately from these formulas. Namely, the group \ $\bar H_8 = \Z$ \ is generated by the cycle $\mbox{
\begin{picture}(21,10)
\put(10.5,3){\oval(21,10)}
\put(3,3){\circle*{2}} 
\put(8,3){\circle*{2}} 
\put(13,3){\circle*{2}} 
\put(18,3){\circle*{2}} 
\end{picture}}_2
+ 
\mbox{
\begin{picture}(21,10)
\put(10.5,3){\oval(21,10)}
\put(3,3){\circle*{2}} 
\put(8,3){\circle*{2}} 
\put(13,3){\circle*{2}} 
\put(18,3){\circle*{2}} 
\end{picture}}_3 + \mbox{
\begin{picture}(21,10)
\put(10.5,3){\oval(21,10)}
\put(3,3){\circle*{2}} 
\put(8,3){\circle*{2}} 
\put(13,3){\circle*{2}} 
\put(18,3){\circle*{2}} 
\end{picture}}_4$ , \ the group \
$H_7 = \Z^2$ \ by cycles 
\begin{equation}
\label{six}
 \mbox{\begin{picture}(16,10)
\put(8,3){\oval(16,10)}
\put(3,1){\circle*{2}} 
\put(3,5){\circle*{2}} 
\put(8,3){\circle*{2}} 
\put(13,3){\circle*{2}} 
\end{picture}}_2+ \mbox{\begin{picture}(16,10)
\put(8,3){\oval(16,10)}
\put(3,1){\circle*{2}} 
\put(3,5){\circle*{2}} 
\put(8,3){\circle*{2}} 
\put(13,3){\circle*{2}} 
\end{picture}}_3 + \mbox{\begin{picture}(16,10)
\put(8,3){\oval(16,10)}
\put(3,1){\circle*{2}} 
\put(3,5){\circle*{2}} 
\put(8,3){\circle*{2}} 
\put(13,3){\circle*{2}} 
\end{picture}}_4- \mbox{\begin{picture}(16,10)
\put(8,3){\oval(16,10)}
\put(3,3){\circle*{2}} 
\put(8,1){\circle*{2}} 
\put(8,5){\circle*{2}} 
\put(13,3){\circle*{2}} 
\end{picture}}_2- \mbox{\begin{picture}(16,10)
\put(8,3){\oval(16,10)}
\put(3,3){\circle*{2}} 
\put(8,1){\circle*{2}} 
\put(8,5){\circle*{2}} 
\put(13,3){\circle*{2}} 
\end{picture}}_3- \mbox{\begin{picture}(16,10)
\put(8,3){\oval(16,10)}
\put(3,3){\circle*{2}} 
\put(8,1){\circle*{2}} 
\put(8,5){\circle*{2}} 
\put(13,3){\circle*{2}} 
\end{picture}}_4+\mbox{\begin{picture}(16,10)
\put(8,3){\oval(16,10)}
\put(3,3){\circle*{2}} 
\put(8,3){\circle*{2}} 
\put(13,1){\circle*{2}} 
\put(13,5){\circle*{2}} 
\end{picture}}_2+\mbox{\begin{picture}(16,10)
\put(8,3){\oval(16,10)}
\put(3,3){\circle*{2}} 
\put(8,3){\circle*{2}} 
\put(13,1){\circle*{2}} 
\put(13,5){\circle*{2}} 
\end{picture}}_3+\mbox{\begin{picture}(16,10)
\put(8,3){\oval(16,10)}
\put(3,3){\circle*{2}} 
\put(8,3){\circle*{2}} 
\put(13,1){\circle*{2}} 
\put(13,5){\circle*{2}} 
\end{picture}}_4 
\end{equation}
and 
\begin{equation}
\label{six2}
 \mbox{\begin{picture}(16,10)
\put(8,3){\oval(16,10)}
\put(3,1){\circle*{2}} 
\put(3,5){\circle*{2}} 
\put(8,3){\circle*{2}} 
\put(13,3){\circle*{2}} 
\end{picture}}_3- \mbox{\begin{picture}(16,10)
\put(8,3){\oval(16,10)}
\put(3,3){\circle*{2}} 
\put(8,1){\circle*{2}} 
\put(8,5){\circle*{2}} 
\put(13,3){\circle*{2}} 
\end{picture}}_3+\mbox{\begin{picture}(16,10)
\put(8,3){\oval(16,10)}
\put(3,3){\circle*{2}} 
\put(8,3){\circle*{2}} 
\put(13,1){\circle*{2}} 
\put(13,5){\circle*{2}} 
\end{picture}}_3 \ ;
\end{equation}
a free generator of the group \ $\bar H_6 = \Z \oplus \Z_2 $ \ can be realized by the cycle \ $\ \mbox{\begin{picture}(12,10)
\put(6,3){\oval(12,10)}
\put(4,0){\circle*{2}} 
\put(4,3){\circle*{2}} 
\put(4,6){\circle*{2}} 
\put(9,3){\circle*{2}} 
\end{picture}}_2 - \mbox{\begin{picture}(12,10)
\put(6,3){\oval(12,10)}
\put(3,3){\circle*{2}} 
\put(8,0){\circle*{2}} 
\put(8,3){\circle*{2}} 
\put(8,6){\circle*{2}} 
\end{picture}}_2$ \ or \ $\ \mbox{\begin{picture}(12,10)
\put(6,3){\oval(12,10)}
\put(4,0){\circle*{2}} 
\put(4,3){\circle*{2}} 
\put(4,6){\circle*{2}} 
\put(9,3){\circle*{2}} 
\end{picture}}_4 \ - \ \mbox{\begin{picture}(12,10)
\put(6,3){\oval(12,10)}
\put(3,3){\circle*{2}} 
\put(8,0){\circle*{2}} 
\put(8,3){\circle*{2}} 
\put(8,6){\circle*{2}} 
\end{picture}}_4$, \ and its element of order 2 by \ $\ \mbox{\begin{picture}(12,10)
\put(6,3){\oval(12,10)}
\put(4,0){\circle*{2}} 
\put(4,3){\circle*{2}} 
\put(4,6){\circle*{2}} 
\put(9,3){\circle*{2}} 
\end{picture}}_3 \ \sim \ \mbox{\begin{picture}(12,10)
\put(6,3){\oval(12,10)}
\put(3,3){\circle*{2}} 
\put(8,0){\circle*{2}} 
\put(8,3){\circle*{2}} 
\put(8,6){\circle*{2}} 
\end{picture}}_3$; \ finally the group \
$\bar H_5 = (\Z_2)^2$ \ is generated by arbitrary two of the three cycles \ $ \mbox{\begin{picture}(9,12)
\put(4,4){\oval(8,13)}
\put(4,-1){\circle*{2}} 
\put(4,2){\circle*{2}} 
\put(4,5){\circle*{2}} 
\put(4,8){\circle*{2}} 
\end{picture}}_2$ \ , \
$ \mbox{\begin{picture}(9,12)
\put(4,4){\oval(8,13)}
\put(4,-1){\circle*{2}} 
\put(4,2){\circle*{2}} 
\put(4,5){\circle*{2}} 
\put(4,8){\circle*{2}} 
\end{picture}}_3$
\ and \ $ \mbox{\begin{picture}(9,12)
\put(4,4){\oval(8,13)}
\put(4,-1){\circle*{2}} 
\put(4,2){\circle*{2}} 
\put(4,5){\circle*{2}} 
\put(4,8){\circle*{2}} 
\end{picture}}_4$\ .

It is easy to calculate (see also \cite{vai}) that groups \ $\bar H_i(B(\C^1,4),\Z)$ \ of the second row of Table \ref{tab5} are generated by the following cycles: \ $\mbox{
\begin{picture}(21,10)
\put(10.5,3){\oval(21,10)}
\put(3,3){\circle*{2}} 
\put(8,3){\circle*{2}} 
\put(13,3){\circle*{2}} 
\put(18,3){\circle*{2}} 
\end{picture}}$ \ for $i=8$, \ $ \mbox{\begin{picture}(16,10)
\put(8,3){\oval(16,10)}
\put(3,1){\circle*{2}} 
\put(3,5){\circle*{2}} 
\put(8,3){\circle*{2}} 
\put(13,3){\circle*{2}} 
\end{picture}} \ - \ \mbox{\begin{picture}(16,10)
\put(8,3){\oval(16,10)}
\put(3,3){\circle*{2}} 
\put(8,1){\circle*{2}} 
\put(8,5){\circle*{2}} 
\put(13,3){\circle*{2}} 
\end{picture}} \ + \ \mbox{\begin{picture}(16,10)
\put(8,3){\oval(16,10)}
\put(3,3){\circle*{2}} 
\put(8,3){\circle*{2}} 
\put(13,1){\circle*{2}} 
\put(13,5){\circle*{2}} 
\end{picture}}$ \ for \ $i=7$, and \ $\mbox{\begin{picture}(8,12)
\put(4,4){\oval(8,13)}
\put(4,-1){\circle*{2}} 
\put(4,2){\circle*{2}} 
\put(4,5){\circle*{2}} 
\put(4,8){\circle*{2}} 
\end{picture}}$ \ for \ $i=5$.

Consider now the exact sequence of all homology groups studied in Theorem \ref{homcal}, defined by the short exact sequence of coefficients \ $\Z \to \widehat{\{A_2\}} \to \{A_2\} $. \ Its map \ $\bar H_i \left(B(\C^1,4), \Z\right) \to \bar H_i (B(\C^1,4), \widehat{\{A_2\}} )$ \ is monomorphic for any \ $i$. \ Namely, it is an isomorphism for \ $i=8$, \ its image for \ $i=7$ \ is generated by cycle (\ref{six}), and for \ $ i=5$ \ by cycle \ $\mbox{\begin{picture}(9,12)
\put(4,4){\oval(8,13)}
\put(4,-1){\circle*{2}} 
\put(4,2){\circle*{2}} 
\put(4,5){\circle*{2}} 
\put(4,8){\circle*{2}} 
\end{picture}}_2 \ - \ \mbox{\begin{picture}(9,12)
\put(4,4){\oval(8,13)}
\put(4,-1){\circle*{2}} 
\put(4,2){\circle*{2}} 
\put(4,5){\circle*{2}} 
\put(4,8){\circle*{2}} 
\end{picture}}_3 \ + \ \mbox{\begin{picture}(9,12)
\put(4,4){\oval(8,13)}
\put(4,-1){\circle*{2}} 
\put(4,2){\circle*{2}} 
\put(4,5){\circle*{2}} 
\put(4,8){\circle*{2}} 
\end{picture}}_4$. This implies the statement of Theorem \ref{homcal} on the structure of groups \ $\bar H_i\left(B(\C^1,4), \{A_2\}\right)$. \hfill $\Box$

\subsection{Homomorphism \ $\partial^1: \E^1_{2,8} \to \E^1_{1,8}$}

By Lemma \ref{lem13}, the source group \ $\E^1_{2,8}$ \ of this homomorphism is naturally isomorphic to the group \ $\bar H_6(B(\C^1,4), \{A_2\})$. \ According to the previous calculation, we can take the cycle \ $\ \mbox{\begin{picture}(12,10)
\put(6,3){\oval(12,10)}
\put(4,0){\circle*{2}} 
\put(4,3){\circle*{2}} 
\put(4,6){\circle*{2}} 
\put(9,3){\circle*{2}} 
\end{picture}}_2 \ - \ \mbox{\begin{picture}(12,10)
\put(6,3){\oval(12,10)}
\put(3,3){\circle*{2}} 
\put(8,0){\circle*{2}} 
\put(8,3){\circle*{2}} 
\put(8,6){\circle*{2}} 
\end{picture}}_2$ \ for its free generator. 

By (\ref{thom3}), the target group \ $\E^1_{1,8}$ \ of this homomorphism is isomorphic to the group \ 
$\bar H_{6}((B(\C^1,3) \times \X) \setminus \triangleq \ , \pm \Z)$. \ By 
exact sequence of the pair \ $(B(\C^1,3) \times \X, \ \triangleq)$ \ (see \S \ref{term31}), any element of this group is characterized by the class of its boundary in the group \ $\bar H_5( \ \triangleq \ , \pm \Z)$. \ By Lemma \ref{lem4}, the last class modulo torsion is characterized by its intersection index in the manifold \ $\triangleq$ \ with the one-dimensional cycle described in statement 5) of this lemma. Let us calculate this intersection index for the boundary of the cycle \ $\partial^1 \left(\left\{ \mbox{\begin{picture}(12,10)
\put(6,3){\oval(12,10)}
\put(4,0){\circle*{2}} 
\put(4,3){\circle*{2}} 
\put(4,6){\circle*{2}} 
\put(9,3){\circle*{2}} 
\end{picture}}_2 - \mbox{\begin{picture}(12,10)
\put(6,3){\oval(12,10)}
\put(3,3){\circle*{2}} 
\put(8,0){\circle*{2}} 
\put(8,3){\circle*{2}} 
\put(8,6){\circle*{2}} 
\end{picture}}_2\right\}\right)$.

\unitlength 0.8mm
\begin{figure}
\begin{picture}(12, 25)
\put(2,0){\circle*{2}}
\put(2,25){\circle*{2}}
\put(2,10){\circle*{2}}
\put(12,15){\circle*{2}}
\put(4,-1){$\gamma$}
\put(4,24){$\alpha$}
\put(4,9){$\beta$}
\put(14,14){$\delta$}
\end{picture} \hspace{6cm}
\begin{picture}(30, 25)
\bezier{350}(2,0)(-4,12.5)(2,25)
\put(2,0){\circle*{2}}
\put(2,25){\circle*{2}}
\put(2,10){\circle*{2}}
\put(12,15){\circle*{2}}
\put(2,0){\line(0,1){10}}
\put(2,0){\line(2,3){10}}
\put(2,10){\line(2,1){10}}
\put(2,25){\line(1,-1){10}}
\put(19,12.5){{\large $-$}}
\end{picture} 
\begin{picture}(12, 25)
\put(2,0){\circle*{2}}
\put(2,25){\circle*{2}}
\put(2,10){\circle*{2}}
\put(12,15){\circle*{2}}
\put(2,0){\line(0,1){10}}
\put(2,25){\line(0,-1){15}}
\put(2,10){\line(2,1){10}}
\put(2,25){\line(1,-1){10}}
\bezier{350}(2,0)(-4,12.5)(2,25)
\end{picture}
\caption{Chain $\ \mbox{ \unitlength 0.5mm \begin{picture}(13,10)
\put(6,3){\oval(13,10)}
\put(4,0){\circle*{2}} 
\put(4,3){\circle*{2}} 
\put(4,6){\circle*{2}} 
\put(9,3){\circle*{2}} 
\end{picture}}_2$}
\label{al3}
\end{figure}

Consider first the contribution of the boundary of the chain \ $\left\{ \mbox{\unitlength 0.5mm \begin{picture}(12,10)
\put(6,3){\oval(12,10)}
\put(4,0){\circle*{2}} 
\put(4,3){\circle*{2}} 
\put(4,6){\circle*{2}} 
\put(9,3){\circle*{2}} 
\end{picture}}_2\right\}$. \
The cell $\ \mbox{ \unitlength 0.5mm \begin{picture}(12,10)
\put(6,3){\oval(12,10)}
\put(4,0){\circle*{2}} 
\put(4,3){\circle*{2}} 
\put(4,6){\circle*{2}} 
\put(9,3){\circle*{2}} 
\end{picture}} \ 
\subset B(\C^1,4)$ \
consists of configurations as in Fig. \ref{al3} (left), i.e. with coinciding three leftmost real values of the corresponding four complex numbers. Recall that we denote these numbers by \ $\alpha, \beta, \gamma, \delta$ \ in such a way that \ $\mbox{Re } \alpha = \mbox{Re } \beta = \mbox{Re } \gamma < \mbox{Re } \delta$, $\mbox{Im } \alpha > \mbox{Im } \beta >\mbox{Im } \gamma$. \ Subscript \ $2$ \ in \ $ \mbox{\unitlength 0.5mm \begin{picture}(12,10)
\put(6,3){\oval(12,10)}
\put(4,0){\circle*{2}} 
\put(4,3){\circle*{2}} 
\put(4,6){\circle*{2}} 
\put(9,3){\circle*{2}} 
\end{picture}}_2$
 indicates the cell over $\ \mbox{ \unitlength 0.5mm \begin{picture}(12,10)
\put(6,3){\oval(12,10)}
\put(4,0){\circle*{2}} 
\put(4,3){\circle*{2}} 
\put(4,6){\circle*{2}} 
\put(9,3){\circle*{2}} 
\end{picture}}$ \ 
in the 3-fold covering of \ $B(\C^1,4)$ \ characterized by matching \ $(\alpha, \beta)(\gamma , \delta)$ \ of these points. This cell corresponds to a $10$-dimensional chain in \ $\Theta_2 \setminus \Theta_1,$ \ which is a fiber bundle over our cell \ $\mbox{ \unitlength 0.5mm \begin{picture}(13,10)
\put(6,3){\oval(12,10)}
\put(4,0){\circle*{2}} 
\put(4,3){\circle*{2}} 
\put(4,6){\circle*{2}} 
\put(9,3){\circle*{2}} 
\end{picture}}$. \ Its fiber over the configuration shown in Fig. \ref{al3} (left) is a chain in the complex of graphs with vertices \ $\alpha, \beta, \gamma, \delta$ \ (see Lemma \ref{lem22}), namely 
the difference of two 4-dimensional simplices represented by graphs obtained from the complete graph by removing edges \ $[\alpha, \beta]$ \ and \ $[\gamma, \delta]$, \ see Fig. \ref{al3} (right).
The boundary of this chain in \ $\Theta_1 \setminus \Theta_0$ \ is the algebraic sum of eight chains corresponding to the three-dimensional boundary faces of these four-dimensional simplices, depicted by not two-connected graphs, see top row of Fig. \ref{fi44}.

\unitlength 0.6mm
\begin{figure}
\mbox{
\begin{picture}(23, 25)
\put(2,0){\circle*{2}}
\put(2,25){\circle*{2}}
\put(2,10){\circle*{2}}
\put(12,15){\circle*{2}}
\put(2,0){\line(0,1){10}}
\put(2,0){\line(2,3){10}}
\put(2,10){\line(2,1){10}}
\put(2,25){\line(1,-1){10}}
\end{picture}
\begin{picture}(23, 25)
\put(2,0){\circle*{2}}
\put(2,25){\circle*{2}}
\put(2,10){\circle*{2}}
\put(12,15){\circle*{2}}
\put(2,0){\line(0,1){10}}
\put(2,0){\line(2,3){10}}
\put(2,10){\line(2,1){10}}
\bezier{350}(2,0)(-4,12.5)(2,25)
\end{picture}
\begin{picture}(23, 25)
\put(2,0){\circle*{2}}
\put(2,25){\circle*{2}}
\put(2,10){\circle*{2}}
\put(12,15){\circle*{2}}
\put(2,0){\line(2,3){10}}
\put(2,10){\line(2,1){10}}
\put(2,25){\line(1,-1){10}}
\bezier{350}(2,0)(-4,12.5)(2,25)
\end{picture}
\begin{picture}(23, 25)
\put(2,0){\circle*{2}}
\put(2,25){\circle*{2}}
\put(2,10){\circle*{2}}
\put(12,15){\circle*{2}}
\put(2,0){\line(0,1){10}}
\put(2,25){\line(1,-1){10}}
\put(2,0){\line(2,3){10}}
\bezier{350}(2,0)(-4,12.5)(2,25)
\end{picture} 
\begin{picture}(23, 25)
\put(2,0){\circle*{2}}
\put(2,25){\circle*{2}}
\put(2,10){\circle*{2}}
\put(12,15){\circle*{2}}
\put(2,25){\line(0,-1){15}}
\put(2,10){\line(2,1){10}}
\put(2,25){\line(1,-1){10}}
\put(2,0){\line(0,1){10}}
\end{picture}
\begin{picture}(23, 25)
\put(2,0){\circle*{2}}
\put(2,25){\circle*{2}}
\put(2,10){\circle*{2}}
\put(12,15){\circle*{2}}
\put(2,0){\line(0,1){10}}
\put(2,25){\line(0,-1){15}}
\put(2,10){\line(2,1){10}}
\bezier{350}(2,0)(-4,12.5)(2,25)
\end{picture}
\begin{picture}(23, 25)
\put(2,0){\circle*{2}}
\put(2,25){\circle*{2}}
\put(2,10){\circle*{2}}
\put(12,15){\circle*{2}}
\put(2,25){\line(1,-1){10}}
\put(2,25){\line(0,-1){15}}
\put(2,25){\line(1,-1){10}}
\put(2,0){\line(0,1){10}}
\bezier{350}(2,0)(-4,12.5)(2,25)
\end{picture}
\begin{picture}(12, 25)
\put(2,0){\circle*{2}}
\put(2,25){\circle*{2}}
\put(2,10){\circle*{2}}
\put(12,15){\circle*{2}}
\put(2,25){\line(1,-1){10}}
\put(2,25){\line(0,-1){15}}
\put(2,10){\line(2,1){10}}
\bezier{350}(2,0)(-4,12.5)(2,25)
\end{picture}
}

\vspace*{0.5cm}
\mbox{
\begin{picture}(21, 40)
\put(2,0){\circle*{2}}
\put(2,11.1){\circle*{2}}
\put(2,10){\circle*{2}}
\put(12,16.1){\circle*{2}}
\put(12,15){\circle*{2}}
\put(2,0){\line(0,1){10}}
\put(2,0){\line(2,3){10}}
\put(2,10){\line(2,1){10}}
\put(2,11.1){\line(2,1){10}}
\put(6,40){\vector(0,-1){10}}
\end{picture}
\begin{picture}(15, 40)
\put(1.2,0){\circle*{2}}
\put(2,0){\circle*{2}}
\put(1.2,10){\circle*{2}}
\put(2,10){\circle*{2}}
\put(12,15){\circle*{2}}
\put(2,0){\line(0,1){10}}
\put(2,0){\line(2,3){10}}
\put(2,10){\line(2,1){10}}
\put(1.2,0){\line(0,1){10}}
\put(6,40){\vector(0,-1){10}}
\end{picture}
\begin{picture}(15, 40)
\put(2,0){\circle*{2}}
\put(2,25){\circle*{2}}
\put(2,-1.5){\circle*{2}}
\put(12,13.5){\circle*{2}}
\put(12,15){\circle*{2}}
\put(2,0){\line(2,3){10}}
\put(2,-1.5){\line(2,3){10}}
\put(2,25){\line(1,-1){10}}
\put(2,0){\line(0,1){25}}
\put(8.5,40){\vector(-1,-4){2.5}}
\end{picture}
\begin{picture}(17, 40)
\put(2,0){\circle*{2}}
\put(2,25){\circle*{2}}
\put(2,26.3){\circle*{2}}
\put(12,15){\circle*{2}}
\put(12,16.3){\circle*{2}}
\put(2,0){\line(2,3){10}}
\put(2,26.3){\line(1,-1){10}}
\put(2,25){\line(1,-1){10}}
\put(2,0){\line(0,1){25}}
\put(3.5,40){\vector(1,-4){2.5}}
\end{picture}
\begin{picture}(23, 40)
\put(2,0){\circle*{2}}
\put(2,25){\circle*{2}}
\put(1,25){\circle*{2}}
\put(1,0){\circle*{2}}
\put(12,15){\circle*{2}}
\put(2,0){\line(0,1){10}}
\put(2,25){\line(1,-1){10}}
\put(2,0){\line(2,3){10}}
\put(1,0){\line(0,1){25}}
\put(2,0){\line(0,1){25}}
\put(4,40){\vector(0,-1){10}}
\end{picture} \hspace{4.3cm}
\begin{picture}(12, 40)
\put(2,25){\circle*{2}}
\put(1,25){\circle*{2}}
\put(2,10){\circle*{2}}
\put(1,10){\circle*{2}}
\put(12,15){\circle*{2}}
\put(2,25){\line(1,-1){10}}
\put(2,25){\line(0,-1){15}}
\put(2,10){\line(2,1){10}}
\put(1,10){\line(0,1){15}}
\put(6,40){\vector(0,-1){10}}
\end{picture}
}
\caption{Boundary of cell \ $\ \mbox{\unitlength 0.5mm \begin{picture}(13,10)
\put(6,3){\oval(13,10)}
\put(4,0){\circle*{2}} 
\put(4,3){\circle*{2}} 
\put(4,6){\circle*{2}} 
\put(9,3){\circle*{2}} 
\end{picture}}_2$ \ in \ $\Theta_1 \setminus \Theta_0$ \ and its boundary over \ $\triangleq$}
\label{fi44}
\end{figure}

Each of these chains is represented by a map to \ $\Theta_1 \setminus \Theta_0$ \ of the space of some fiber bundle over the cell 
\unitlength 0.5mm
$\ \mbox{\begin{picture}(12,10)
\put(6,3){\oval(12,10)}
\put(4,0){\circle*{2}} 
\put(4,3){\circle*{2}} 
\put(4,6){\circle*{2}} 
\put(9,3){\circle*{2}} 
\end{picture}} $, 
whose fiber over the configuration \ $\Gamma=\{\alpha, \beta, \gamma, \delta\} \in \ 
\mbox{\begin{picture}(12,10)
\put(6,3){\oval(12,10)}
\put(4,0){\circle*{2}} 
\put(4,3){\circle*{2}} 
\put(4,6){\circle*{2}} 
\put(9,3){\circle*{2}} 
\end{picture}}$
\ is the corresponding three-dimensional face of the 5-simplex \ $\Delta(\Gamma)$. \ This map takes this face isomorphically to a three-dimensional simplex which is the fiber of the bundle \ $\Theta_1 \setminus \Theta_0 \to (B(\C^1,3) \times \X ) \setminus \triangleq$ \ over a certain point of this base space, depending on both the configuration \ $\{\alpha, \beta, \gamma, \delta\}$ \ and the mapped face. For example, the first picture in the top row of Fig. \ref{fi44} consists of a triangle with vertices \ $\beta, \gamma, \delta$ \ and the additional edge \ $(\alpha, \delta)$, \ therefore the corresponding face of the 5-simplex over the configuration \ $\{\alpha, \beta, \gamma, \delta\} \in B(\C^1,4)$ \ goes to the fiber over the point \ $\{\beta, \gamma, \delta\} \times \{\alpha, \delta\} \in B(\C^1,3) \times \X$. 

Boundary points in \ $\triangleq$ \ of either of these eight chains are approached in the degeneration of underlying configurations \ \unitlength 0.5mm $\{\alpha, \beta, \gamma, \delta\} \in \ \mbox{\begin{picture}(12,10)
\put(6,3){\oval(12,10)}
\put(4,0){\circle*{2}} 
\put(4,3){\circle*{2}} 
\put(4,6){\circle*{2}} 
\put(9,3){\circle*{2}} 
\end{picture}} $, \ when the only 1-valent vertex of the corresponding graph tends to some other vertex (not joined with it by an edge). All these degenerations are shown in the bottom row of Fig. \ref{fi44}.

 The limit configuration obtained by such a degeneration can belong to the basic 1-cycle of the group \ $H_1(\triangleq,\pm \Z)$ \ indicated in statement 5) of Lemma \ref{lem4}
only if it is the configuration \unitlength 0.8 mm \
\begin{picture}(10,7)
\put(0,5){\circle*{1.5}}
\put(0,0){\circle*{1.5}}
\put(10,5){\circle*{1.5}}
\put(0,6){\circle*{1.5}}
\put(10,6){\circle*{1.5}}
\put(0,0){\line(0,1){5}}
\put(0,5){\line(1,0){10}}
\put(0,0){\line(2,1){10}}
\put(0,6){\line(1,0){10}}
\end{picture} \ \
or \ \
\begin{picture}(10,7)
\put(0,1){\circle*{1.5}}
\put(0,0){\circle*{1.5}}
\put(10,1){\circle*{1.5}}
\put(0,6){\circle*{1.5}}
\put(10,0){\circle*{1.5}}
\put(0,1){\line(0,1){5}}
\put(0,1){\line(1,0){10}}
\put(0,6){\line(2,-1){10}}
\put(0,0){\line(1,0){10}}
\end{picture} \ \ with vertices of triangles at the points $0$ (right angle), 1 and $\mp i/2$. These points appear only in three pieces shown in the bottom row of Fig. \ref{fi44}: the first, third and fourth; each of them defines a transversal intersection of the corresponding cycles and makes contribution \ $\pm 1$ \ to the intersection index. 

Similar considerations with the summand \ \unitlength 0.5mm $\ \mbox{\begin{picture}(12,10)
\put(6,3){\oval(12,10)}
\put(3,3){\circle*{2}} 
\put(8,0){\circle*{2}} 
\put(8,3){\circle*{2}} 
\put(8,6){\circle*{2}} 
\end{picture}}_2$ \ of our generator of the group \ $\E^1_{2,8}$ \ give us no boundary components in \ $\Theta_1$ \ that might contribute to the intersection index with this 1-cycle. Thus, this index is an odd number, and the last statement of Theorem \ref{thmf3} is proved. \hfill $\Box$

\section{Estimates in the stable spectral sequence}
\label{estim}

\begin{proposition}
\label{pro99}
For any natural \ $\rho$ \ in stable range $($i.e. satisfying $($\ref{stabdim}$))$
and any \ $j \in \{0, 1, \dots, \rho-1\}$, 
 the Borel--Moore homology group $$\bar H_i(\Theta_j \setminus \Theta_{j-1})$$ of term \ $\Theta_j \setminus \Theta_{j-1}$ \ of the inverse filtration of term $\Lambda_\rho \setminus \Lambda_{\rho-1}$ of main filtration of \ $\Lambda$ \
$($or, which is the same, group \ $\E^1_{j, i-j}$ \ of the auxiliary spectral sequence defined by this inverse filtration$)$ is trivial for all \ $i > 5\rho -j-1 $.
\end{proposition}

\noindent
{\it Proof.} By Lemma \ref{lem18} and Corollary \ref{cor1}, \ $\Theta_j \setminus \Theta_{j-1}$ \ is the union of spaces associated with sets \ $A$ \ such that \ $|A|-\#(A) = \rho$ \ and \ $j = |A|-2\#(A)$. \ Since \ $\rho$ \ is in stable range, each of these spaces is a fiber bundle, whose base is \ $2|A|$-dimensional, and the Borel--Moore homology group of the fiber is non-trivial only in the dimension \ $2|A|-3\#(A) -1$ \ (see Corollary \ref{cor1}). Therefore the Borel--Moore homology group of any such space is trivial in dimensions greater than \ $4|A|-3\#(A)-1 \equiv 5\rho - j - 1$. \hfill $\Box$ 

\begin{lemma}[see e.g. \cite{coh}]
\label{lem80}
If \ $\rho >1$ \ then the group \ $\bar H_i(B(\X,\rho), \pm \Z)\equiv \bar H_i(B(\R^4,\rho), \pm \Z)$ \ is finite for any \ $i$, \ trivial for \ $i=4\rho$, \ and isomorphic to \ $\Z_2$ \ for \ $i=4\rho-1$. \hfill $\Box$
\end{lemma}

\begin{proposition}
\label{pp0}
For any \ $\rho>1$ \ in the stable range, the Borel--Moore homology group \ $\bar H_i(\Theta_0)$ \ of term \ $\Theta_0$ \ of \ $\Lambda_\rho \setminus \Lambda_{\rho-1}$ \ is 

a$)$ trivial for \ $i> 5\rho-2$,

b$)$ isomorphic to \ $\Z_2$ \ for \ $i=5\rho - 2,$ and 

c$)$ finite for \ $i > 5\rho - 6$.
\end{proposition}

\noindent
{\it Proof.} This term contains only one \ $A$-block, where \ $A = (2, 2, \dots , 2)$ ($\rho$ \ deuces). Similarly to \S \ref{t30}, 
$$\bar H_i(\Theta_0) \simeq \bar H_{i- (\rho-1)}(\tilde B(\X,\rho), \pm \Z), $$
where \ $\tilde B(\X,\rho)$ \ is obtained from the configuration space \ $B(\X,\rho)$ \ by removing a subset of complex codimension \ $\geq 3$ \ (consisting of non-regular \ $\bar A$-configurations). All statements of the proposition follow from the homological exact sequence of the pair \ $(B(\X, \rho), B(\X, \rho) \setminus \tilde B(\X, \rho))$, \ Lemma \ref{lem80} and dimensional restrictions on homology groups of the removed subset. \hfill $\Box$ 

\begin{proposition}
\label{pp1}
For any \ $\rho >1$ \ in the stable range, the Borel--Moore homology group of term \ $\Theta_1 \setminus \Theta_0$ \ of \ $\Lambda_\rho \setminus \Lambda_{\rho-1}$ \ is 

a$)$ trivial in dimensions exceeding \ $5\rho -3 $, 

b$)$ finite in dimensions exceeding \ $5\rho - 6$.
\end{proposition}

\noindent
{\it Proof.} The unique \ $A$-block covering this term is \ $\{3, 2 , 2, \dots, 2\}$ \ ($\rho -2$ \ deuces). This term is the space of a fiber bundle with \ $2(2\rho-1)$-dimensional base \ $\tilde \B(A)$ \ and fibers equal to \ $\rho$-dimensional simplices. 
This base is orientable, but the orientation of fibers is violated by some loops in the base, therefore the homology group in top dimension \ $5\rho -2$ \ is trivial, which implies statement a). 

The base space \ $\tilde \B(A)$ \ is the direct product \ $B(\C^1,3) \times B(\X, \rho -2)$, \ from which a subset \ $\yen$ \ of complex codimension \ $2$ \ is removed. The local system on this base formed by \ $\rho$-dimensional homology groups of fibers is isomorphic to one induced from the tensor product of local systems \ $\pm \Z$ \ on the factors \ $B(\C^1,3)$ \ and \ $B(\X,p-2)$. Denote the last tensor product by \ $\gimel$. \ By the K\"unneth formula and Lemma \ref{lem6}, the homology group of the product \ $ B(\C^1,3) \times B(\X, \rho -2)$ \ with these coefficients is finite, so by Thom isomorphism and exact sequence of pair $ ((B(\C^1,3) \times B(\X, \rho -2)), \yen)$
we get isomorphism 
\begin{equation}
\label{rdu5}
\bar H_i(\Theta_1 \setminus \Theta_0, \Q) \simeq \bar H_{i-\rho}(\tilde \B(A), \gimel \otimes \Q)\simeq \bar H_{i-(\rho+1)}(\yen, \gimel \otimes \Q).
\end{equation} 
 The dimension of \ $\yen$ \ is equal to \ $4\rho - 6$, which implies the triviality of groups (\ref{rdu5}) for all \ $i > 5\rho -5$. \ To also overcome the top dimension \ $5 \rho -5$, \ note that a Zariski open subset of \ $\yen$ \ coincides with a Zariski open subset of the direct product of the space \ $\triangleq$ \ considered in \S \ref{term31} and configuration space \ $B(\X, \rho -2)$. \ This open set is orientable, but the restriction to it of the one-dimensional local system \ $\gimel$ \ is not constant, so the group (\ref{rdu5}) with \ $i=2(2\rho -3)+(\rho +1) \equiv 5\rho -5$ \ is also trivial.
 \hfill $\Box$

\begin{proposition}
\label{pp4}
For any \ $\rho >3$ \ in the stable range, 

a$)$ Borel--Moore homology group \ $\bar H_i(\Theta_2 \setminus \Theta_1)$ \ of term \ $\Theta_2 \setminus \Theta_{1}$ \ of inverse filtration of \ $\Lambda_\rho \setminus \Lambda_{\rho-1}$ \ is 
trivial if \ $i >5\rho -4$; 

b$)$ if \ $\rho >4$, \ then this group is finite for any \ $i> 5\rho -6$; 

c$)$ in the case \ $\rho=4$, \ $\bar H_{16}(\Theta_2 \setminus \Theta_1) \simeq \Z$.
\end{proposition}

\noindent
{\it Proof.} This term \ $\Theta_2 \setminus \Theta_1$ \ is covered by two \ $A$-blocks corresponding to \ $A = \{4, 2, 2, \dots, 2\}$ ($\rho -3$ \ deuces) and \ $A=\{3, 3, 2, \dots, 2\}$ ($\rho - 4$ \ deuces). The first \ $A$-block is a fiber bundle over the base \ $\tilde \B(A)$ \ equal to the product $B(\C^1,4) \times B(\X, \rho - 3)$, \ from which a subset of complex codimension \ $2$ \ is removed. By Theorem \ref{2c}, the Borel--Moore homology group of any of its fibers is isomorphic to \ $\Z^2$ \ in dimension \ $\rho +1$ \ and is trivial in all other dimensions. The local system on \ $\tilde B(A)$ \ formed by these homology groups of fibers is isomorphic to the tensor product of the local systems induced from local systems \ $\{A_2\}$ \ and \ $\pm \Z$ \ on the factors of \ $B(\C^1,4) \times B(\X, \rho -3)$. The estimates of Proposition \ref{pp4} applied exclusively to the homology groups of this \ $A$-block follow now from the exact sequence of the pair \ $(B(\C^1,4) \times B(\X, \rho -3), \{\mbox{the removed set}\})$, \ statements of Theorem \ref{homcal} and Lemma \ref{lem80} on homology of factors, the K\"unneth formula, and dimensional restrictions on the homology groups of the removed set.

In the case of the \ $A$-block with \ $A=(3, 3, 2, ..., 2)$, \ the proof is similar, with Theorem \ref{homcal} replaced by Lemma \ref{lem6} and additional accounting of the factorization by the permutation of two sets of cardinality \ 3 \ of \ $\bar A$-configurations (which does not increase the rational homology groups). \hfill $\Box$

\subsection{Proof of statement 4 of Theorem \ref{mthm2}}
\label{prof4}

According to Propositions \ref{pp0} (c), \ref{pp1} (b), \ref{pp4} (a), and \ref{pro99}, the following groups \ $\E^1_{j, \mathfrak {q} } $ \ of the spectral sequence defined by the inverse filtration of the term \ $\Lambda_4 \setminus \Lambda_3$ \ are finite: $\E^1_{0, \mathfrak q}$ \ for \ $ \mathfrak q > 14$, \ $\E^1_{1, \mathfrak q}$ \ for \ $ \mathfrak q > 13$, \ $\E^1_{2, \mathfrak q}$ \ for \ $ \mathfrak q > 14$, \ and \ $\E^1_{j, \mathfrak q}$ \ for \ $ \mathfrak q > 19-2j$ \ and all \ $j$. By Proposition \ref{pp4} (c), group \ $\E^1_{2, 14}$ \ is equal to \ $\Z$. 

Thus all infinite groups \ $\E^1_{j, \mathfrak {q} } $ \ with \ $j+ \mathfrak {q} \geq 16$ \ are only \ $\E^1_{2,14} \simeq \Z$ \ and possibly \ $\E^1_{3,13} $ \ (which is presumably also finite). Both differentials \ $\partial^1: \E^1_{2,14} \to \E^1_{1,14}$ \ and \ $\partial^2: \E^2_{2,14} \to \E^2_{0,15}$ \ act into finite groups. Therefore the group $$\bar H_{16}(\Lambda_4 \setminus \Lambda_3) \simeq \bar H_{2D-8k+16}(\sigma_4 \setminus \sigma_3) \equiv E^1_{4,2D-8k+12}$$ is infinite.
By (\ref{aldu}), this implies statement 4 of Theorem \ref{mthm2}. \hfill $\Box$

\subsection{Proof of statement 5 of Theorem \ref{mthm2}}
\label{prof5}

a) By Propositions \ref{pro99} (applied to all \ $j \geq 1$) and \ref{pp0} (a), all groups \ $\E^1_{j, \mathfrak {q} } $ \ with \ $j+ \mathfrak{q}> 5\rho -2$ \ of the spectral sequence defined by the inverse filtration of any term \ $\Lambda_\rho \setminus \Lambda_{\rho-1}$ \ are trivial. Hence the groups 
\begin{equation}
\label{lasteq}
\bar H_i(\Lambda_\rho \setminus \Lambda_{\rho-1}) \simeq \bar H_{2D-2k \rho +i}(\sigma_\rho \setminus \sigma_{\rho-1}) \equiv E^1_{\rho, 2D-(2k+1)\rho+i}
\end{equation}
 are also trivial for \ $i>5\rho-2$. 

b) By Propositions \ref{pro99} (applied to all \ $j \geq 2$), \ref{pp0} (b) and \ref{pp1} (a), the only such non-trivial group with $j+ \mathfrak {q}=5\rho -2$ \ is \ $\E^1_{0, 5\rho -2 } \simeq \Z_2$. \ Therefore the group $$\bar H_{5\rho-2}(\Lambda_\rho \setminus \Lambda_{\rho-1}) \simeq \bar H_{2D-(2k-5)\rho-2}(\sigma_\rho \setminus \sigma_{\rho-1}) \equiv E^1_{\rho, 2D-(2k-5)\rho-2} $$
is isomorphic to \ $\Z_2$.

c) By Proposition \ref{pro99} (applied to all \ $j \geq 3$), \ref{pp0} (c), \ref{pp1} (b) and \ref{pp4} (b), all such groups with \ $j + \mathfrak q > 5\rho -4$ \ are finite. Therefore groups (\ref{lasteq}) are also finite for \ $i> 5\rho-4$.
\smallskip

Transformation (\ref{aldu}) turns these three facts into three assertions of statement 5 of Theorem \ref{mthm2}. \hfill $\Box$

\section{Two problems}

1. Find an interpretation and a combinatorial formula for the basic element of stable group \ $H^{6k-12}(P(\infty,k) \setminus \Sigma, \Q)$ \ (at least for \ $k=3$).

2. Calculate the multiplication $$H^{2k-5}(P(\infty,k) \setminus \Sigma, \Q) \otimes H^{6k-12}(P(\infty,k) \setminus \Sigma, \Q) \to H^{8k-17}(P(\infty, k) \setminus \Sigma, \Q).$$


\begin{thebibliography}{99}

\bibitem{AA} V.I.~Arnold, On some topological invariants of algebraic functions, Trans. Moscow Math. Society, 1970, V.21, 30--52.

\bibitem{AVG} 
V.I.~Arnold, A.N.~Varchenko, and S.M.~Gusein-Zade, {\it Singularities of differentiable maps}, vol.1: {\it Classification of critical points, caustics, and wavefronts}. ``Nauka", Moscow, 1982; Engl. transl.: Birkh\"auser, Basel, 2012.

\bibitem{A} V.I.~Arnold a.o., {\it Arnold's problems}. Springer and Phasis, 2004.

\bibitem{BBB} E. Babson, A. Bj\"orner, S. Linusson, J. Shareshian, and V. Welker, {\em Complexes
of not i-connected graphs}, Topology 38:2 (1999), 271–299.

\bibitem{coh} F.R.~Cohen, {\it The homology of ${\mathcal C}_{n+1}$-spaces,} in: F.R.~Cohen, T.J.~Lada, J.P.~May, {\it The homology of iterated loop spaces}, Lect. Notes Math., Berlin, Springer--Verlag, {\bf 533} (1976), pp. 207--353.

\bibitem{fuks} D.B.~Fuchs, {\it Cohomology of braid groups mod 2}, Funct. Anal. Appl., 4:2 (1970), 46--59.

\bibitem{GM} M.~Goresky, R.~MacPherson, {\it Stratified Morse Theory,} Springer, 1988, Berlin a.o. 

\bibitem{KhR} B.~Khesin, A.~Rosly, {\it Polar homology and holomorphic bundles}, 
Philos. Trans. of the Royal Society of London. Ser. A, 359:1784 (2001), 1413--1427.

\bibitem{koz} D.~Kozlov, {\it A comparison of Vassiliev and Ziegler–Živaljević models for homotopy types of subspace arrangements}, 
Topology and its Applications, 126:1--2 (2002), 119--129.


\bibitem{tur} V.E.~Turchin, {\em Homologies of complexes of doubly connected graphs,} Russian
Math. Surveys (Uspekhi) 52 (1997), 426–427.

\bibitem{ks} V.A.~Vassiliev, {\it Cohomology of knot spaces}, in: Theory of Singularities and its Applications, Adv. Soviet Math., 1, AMS, Providence, RI, 1990, 23–69.

\bibitem{book} V.A.~Vassiliev, {\it Complements of Discriminants of Smooth Maps: Topology and Applications}, 2nd edition, Trans. Math. Monogr., 98, AMS, Providence, RI, 1994.


\bibitem{twocon} V.A.~Vassiliev, {\it Topology of two-connected graphs and homology of spaces of knots}, in: Differential and symplectic topology of knots and curves, Amer. Math. Soc. Transl. Ser. 2, 190, AMS, Providence, RI, 1999, 253–286.

\bibitem{vai} F.V.~Weinstein, {\em Cohomologies of braid groups}, Funct. Anal. Appl., 12:2 (1978), 135–137.
\end{thebibliography}
\end{document}